\documentclass[final,onefignum,onetabnum]{siamart190516}

\usepackage{cite}
\usepackage{enumitem}
\usepackage[utf8]{inputenc}
\usepackage[english]{babel}
\usepackage[a4paper]{geometry}

\usepackage[normalem]{ulem}

\usepackage{amssymb}
\usepackage{amsmath}

\newtheorem*{remark}{Remark}
\usepackage[algo2e,vlined,ruled]{algorithm2e}
\usepackage{multirow}
\usepackage{tcolorbox}

\newcommand{\bfx}{\mathbf{x}}
\newcommand{\bfr}{\mathbf{r}}
\newcommand{\bfu}{\mathbf{u}}
\newcommand{\bfg}{\mathbf{g}}
\newcommand{\bfe}{\mathbf{e}}

\newcommand{\bfb}{\mathbf{b}}

\newcommand{\bff}{\mathbf{f}}
\newcommand{\norm}[1]{\left\lVert#1\right\rVert}
\usepackage{graphicx}
\usepackage{subcaption}
\graphicspath{ {./images/} }
\usepackage{wrapfig}
\usepackage[export]{adjustbox}
\usepackage{hyperref}
\numberwithin{equation}{section}

\headers{Deep Learning Preconditioners for the Helmholtz Equation}{Yael Azulay and Eran Treister}

\title{Multigrid-augmented deep learning preconditioners for the Helmholtz equation\thanks{Corresponding author: Eran Treister. \funding{This research was supported by The Israel Science Foundation (grant No. 1589/19), and in part by the Israeli Council for Higher Education (CHE) via the Data Science Research Center, BGU, Israel.}}}

\author{Yael Azulay\thanks{Ben-Gurion University, Beer Sheva, Israel.
  \email{yaelaz@post.bgu.ac.il, erant@cs.bgu.ac.il}}
\and Eran Treister$^\dag$}

\date{}

\begin{document}
\maketitle

\section*{Abstract}
In this paper, we present a data-driven approach to iteratively solve the discrete heterogeneous Helmholtz equation at high wavenumbers. In our approach, we combine classical iterative solvers with convolutional neural networks (CNNs) to form a preconditioner which is applied within a Krylov solver. For the preconditioner, we use a CNN of type U-Net that operates in conjunction with multigrid ingredients. Two types of preconditioners are proposed 1) U-Net as a coarse grid solver, and 2) U-Net as a deflation operator with shifted Laplacian V-cycles. Following our training scheme and data-augmentation, our CNN preconditioner can generalize over residuals and a relatively general set of wave slowness models. On top of that, we also offer an encoder-solver framework where an ``encoder'' network generalizes over the medium and sends context vectors to another ``solver'' network, which generalizes over the right-hand-sides. We show that this option is more robust and efficient than the stand-alone variant. Lastly, we also offer a mini-retraining procedure, to improve the solver after the model is known. This option is beneficial when solving multiple right-hand-sides, like in inverse problems. We demonstrate the efficiency and generalization abilities of our approach on a variety of 2D problems.


\begin{keywords}
  Helmholtz equation, Shifted Laplace preconditioner, Iterative solution methods, Multigrid, Deep Learning, U-Net, Convolutional Neural Networks.
\end{keywords}

\begin{AMS}
  68T07, 65N55, 65N22
\end{AMS}

\section{Introduction}
The efficient numerical solution of the Helmholtz equation with a high and spatially-dependent wavenumber is a difficult task. The linear system that results from the discretization of the Helmholtz equation involves a large, complex-valued, sparse, and indefinite matrix. As the wavenumber is higher, the linear system needs to be larger, and the solution becomes more oscillatory and spatially complex. Many authors have contributed to the development of various iterative methods for solving the problem, e.g.,  
\cite{poulson2013parallel,luo2014fast,gander2019review,graham2020domain,SHEIKH2016473,dwarka2020scalable,treister2019multigrid,erlangga2006novel,olson2010smoothed,10.1145/3054946},
yet it remains expensive and difficult to solve.

One of the most common methods to solve the heterogeneous Helmholtz equation is the Shifted Laplacian (SL) multigrid method \cite{erlangga2006novel,umetani2009multigrid}, which is based on a multigrid preconditioner for a damped version of the equation. That is, a complex shift is added to the matrix (and to its eigenvalues) to form a preconditioner that is efficiently solved by standard multigrid components. The preconditioning is typically applied with a Krylov solver, such as (flexible) GMRES \cite{saad1993flexible}. However, the shifted Laplacian preconditioner is not scalable due to the high shift that is required for multigrid to solve the shifted system, and further investigation on the topic is needed. The work in
\cite{sheikh2013convergence} provides analysis of the SL solver, showing that the increase of the wavenumber leads to a large increase in the number of iterations and the cost required for the solution. Indeed, the number of iterations in SL increases linearly with the wavenumber. 

The numerical methods mentioned above are in some sense ``classical''. In recent years, data-driven methods have revolutionized many fields, and in particular, Deep Learning (DL) approaches have changed the landscape of fields like computer vision, natural language processing and bioinformatics \cite{lecun2015deep,krizhevsky2012imagenet}.
Similarly, solving partial differential equations (PDEs) with a deep neural network is a natural idea, and has been recently considered in various forms, which are roughly divided into two, as we present below. Another related approach is optimization of PDE solvers \cite{greenfeld2019learning}.

One common class of DL-for-PDEs approaches, which are often referred to as physics-informed neural networks (PINNs), includes both PDE solvers  \cite{raissi2018hidden,raissi2019physics,lu2020deepxde, han2018solving},  and PDE discovery methods \cite{bar2019learning}.
The general idea is to implement a neural network that implicitly represents the PDE solution.
That is, the solution is approximated by the neural network as an analytical function at a given point in space. This function is obtained using stochastic optimization for learning the network's weights as an unsupervised learning task, given a large training set of points. The loss function that is minimized forces the network to meet both the system of equations and the initial and/or boundary conditions, using automatic differentiation to get the PDE derivatives. Therefore, the learning process essentially minimizes the residual, which is computed by differentiating the neural network function to get the terms of the PDE analytically.
Following the learning process, the PDE coefficients and boundary conditions are assumed to be embedded into the neural network weights. This approach can be applied to different types of PDEs. As PINNs are evaluated at a single point, their training is in principle mesh-free, but requires the training data set and the architecture to be large enough to capture the complexity of the solution and the problem. A nice property of PINNs is that they can also solve inverse problems, as suggested, for example, in \cite{bar2019unsupervised} for the weighted Poisson equation and its inverse problem, using extended unsupervised loss function that includes the residual in $L_{\infty}$ norm.

A different family of approaches involves solvers in which DL is used to approximate the solution of the PDE in an end-to-end fashion given the problem properties (coefficients, initial conditions or sources). That is, the network needs to generalize over some of the PDE properties, which is the more common scenario in computer vision applications. For example, in \cite{Hermann2020} the authors propose 'PauliNet', a  wavefunction ansatz based on a neural network that achieves nearly exact solutions of the electronic Schrödinger equation. 
Similarly, some works employ Convolutional Neural Networks (CNNs) which are typically used for applications involving data on regular grids, like images and videos, or, discretized PDEs on a structured grid. For example, the authors of \cite{geist2021numerical, khoo2021solving} suggested CNNs to solve 1D and 2D parametric elliptic PDEs. Khoo and Ying \cite{khoo2018switchnet} introduced a CNN named 'SwitchNet', for approximating forward and inverse maps arising from the time-harmonic wave equation, which is the Helmholtz equation in time domain. In these problems, the local information in the input has a global impact on the output, hence, SwitchNet is based on convolution layers and operators of reshaping, switching, flattening, and point-wise matrix multiplication, with which the network performs knowledge-based actions and mimics the butterfly Helmholtz preconditioner \cite{li2015butterfly}. The network is trained to solve the problem for plane-wave sources and generalize for 2-4 Gaussian scatterers in a constant medium. Similarly, \cite{wiecha2019deep} employs a 3D CNN to predict near and far scattering fields. However, no neural network has been presented in the literature that effectively solves the Helmholtz equation at high wavenumber for rather general heterogeneous models and right-hand-sides.

The main principle of our work is merging data-driven approaches, such as DL, together with classical methods. We wish to exploit the strengths of the classical solvers and to complement them with DL, covering for their weaknesses. One DL accelerator adopting this approach, named CFDNet \cite{Obiols_Sales_2020}, has been suggested for fluid simulations. The authors present a coupled framework using both simulation and deep learning, for accelerating the convergence of Reynolds averaged Navier-Stokes simulations. They show that their model performs well and also generalizes for examples not seen in the training. We note that while DL models are considered to be computationally heavy, the emergence of graphics processing units (GPUs) makes DL efficient and beneficial, both in terms of computational cost and in terms of parallelism. Furthermore, compression methods can really accelerate the inference of DL methods by a significant factor \cite{cheng2018model}, using a variety of approaches like pruning (sparsification), quantization, knowledge distillation and more, as the community aims to make DL feasible for low-power edge devices.

In order to solve the Helmholtz equation for a heterogeneous medium and a high wavenumber in a single application of a neural network, one would have to train a huge network to be able capture all the complexity in the solution---a lot of oscillatory waves, reflections and interference. That is even for a certain type of right-hand-sides. In the context of PINNs, the high variability of the solution may require the training to include a huge amount of points, which, together with the fact that we need a large network, makes the training expensive at solve time. To be effective, PINNs assume that the PDE solution is relatively simple, and that the network is light and can be trained with a relatively small data set of samples. As example, the work of \cite{moseley2020solving} considered the time-domain wave equation, and limited the solution to point sources, and a given heterogeneous medium (i.e., the network does not generalize over the medium---it is fixed). While the first arrival wave is indeed well predicted, it seems that reflections are not recovered well by the NN, compared to the true solution. These reflections are important to capture for geophysical applications, e.g., Full Waveform Inversion \cite{metivier2017review}. Lastly, we note that NNs can easily approximate simple functions or solve problems to reasonable-but-low accuracy using a relatively low number of parameters. However, when high complexity and accuracy are needed, the typical size of the neural networks that are required grow tremendously. For example, to improve classification accuracy by 2\% only, architectures need to almost double their size \cite{he2015deep}.

To summarize, given a problem we cannot be sure that the architecture is sufficient to accurately represent the solution, and in hard cases, huge architectures are needed. This brings us to the conclusion that in order to accurately solve complicated problems, the network has to be used iteratively to build and improve the solution in several steps, like a preconditioner. This requires, however, the ability to work on general right-hand-sides (residuals) which are multivariate grid-shaped inputs.  Furthermore, if we wish to  \textit{generalize} on the medium as an input of the network, we have to use a network that consumes images---a CNN---as a full heterogeneous medium is an image. We note that generating training data for our task and training the CNN model are very expensive, and we assume they are done once off-line to be used multiple times later. However, the preconditioner is fast and efficient in the solution time once the model is trained. That is also a reason it is important to have the slowness model generalization.

For the reasons stated above, in this work we utilize CNNs as preconditioners for the discretized Helmholtz equation. We propose to use a CNN together with a classical multigrid approach---the Shifted Laplacian---which is efficient at removing only part of the error. We build on the similarity between a V-cycle and a U-Net \cite{ronneberger2015unet} which is a multiscale CNN used for image-to-image mappings such as semantic segmentation, image denoising etc. Indeed, a geometric V-cycle can be naturally implemented using DL frameworks, suggesting that a U-Net architecture, like a V-cycle, can be used as a preconditioner for solving PDEs on a regular grid. Also, a V-cycle can be used to solve the equation for any medium and right-hand-side, meaning that it generalizes on the problem. This suggests that a U-Net or at least a variant of it can be trained to generalize on both a right-hand-size and a medium, and act as a preconditioner. In this work we consider a standard U-Net, but note that other variants are worthy of consideration.


To precondition the Helmholtz equation, we use a U-Net together with a classical approach, to complement the network and stabilize the solution process. Our U-Net architecture begins with a strided convolution, initially projecting the fine level input to a coarse grid. We examine two preconditioners: one is a U-Net that acts as a non-linear coarse grid correction and is applied with a Jacobi iteration for pre-and post-smoothing.
In the second formulation we apply alternating U-Net and SL multigrid V-cycle, which provides a more powerful smoothing and somewhat resembles the deflation preconditioner of \cite{sheikh2013convergence,dwarka2020scalable}. Both preconditioners are applied within a Krylov subspace method---flexible GMRES (FGMRES). We define our U-Net to generalize over a random right-hand-side and a random piece-wise smooth slowness model. We also suggest two upgrades to the scheme above: (1) an encoder-solver framework and (2) a mini-retrain approach. Both upgrades are detailed later, and can improve the standard U-Net tremendously.




\section{Preliminaries and Background}

\subsection{Problem formulation}

The heterogeneous Helmholtz problem is give by
\begin{equation}\label{eq:helmholtz}
  -\Delta u(\vec x,\omega)-\omega^2\kappa(\vec x)^2(1-\gamma i)u(\vec x,\omega)=g(\vec x,\omega), \quad  \vec x\in \Omega
\end{equation}
The unknown $u(\vec x,\omega)$ represents the pressure wave function in the frequency domain, $\omega=2\pi f$ denotes the angular frequency,
$\Delta$ is the Laplacian operator and $\kappa(\vec x)$  is the heterogeneous wave slowness model - the inverse of its velocity.
The source term is the right-hand-side (RHS) of the system, and is denoted by $g(\vec x,\omega)$. The parameter $\gamma$ indicates the fraction of attenuation (or damping) in the medium and $i =\sqrt{-1}$ . We focus here on the hardest case for iterative methods, which is $\gamma=0$, but it is also possible to consider a varying attenuation. As boundary conditions, we consider an absorbing layer \cite{erlangga2006novel,engquist1979radiation}, but Sommerfeld, PML \cite{singer2004perfectly} or \cite{papadimitropoulos2021double} can be viable options as well.

\subsubsection*{Discretization}
The second order finite difference discretization of the aforementioned problem \eqref{eq:helmholtz} on a uniform mesh of width $h$ in both $x$ and $y$ directions yields the stencil
\begin{equation}\label{eq:ah_matrix}
 A^h=\frac{1}{h^2}
 \begin{bmatrix}
 0 & -1 & 0\\
 -1 & 4-\omega^2 \kappa(\mathbf{x})^2h^2 & -1\\
 0 & -1 & 0
 \end{bmatrix},
\end{equation}
where boldface letters like $\bfx$ denote discrete vectors. This leads to a system of linear equations
\begin{equation}\label{eq:aug}
A^h \bfu^h =\bfg^h. \end{equation}
When discretizing the Helmholtz equation, we must obey the rule of thumb that at least 10 grid nodes per wavelength are used, which leads to the bound $\omega \kappa h \leq \frac{2\pi}{10} \approx 0.628$. This typically requires a very fine mesh for high wavenumbers, significantly increasing the number of unknowns and making the system  \eqref{eq:aug} huge, ill-conditioned, in addition to being indefinite, and complex-valued due to the boundary conditions. Moreover, the number of eigenvalues of the matrix $A^h$ in \eqref{eq:ah_matrix} with negative real part increases as the wavenumber $\kappa\omega$ is higher.
Hence, solutions for two and especially three dimensional problems are challenging and require creative iterative solvers.


\subsection{Shifted Laplacian multigrid}

The common approach to solving \eqref{eq:aug} is to use preconditioner within an iterative Krylov subspace method.  
The SL operator
\begin{equation}\label{eq:mu_sl}
 Mu = -\Delta u-\omega^2\kappa(\vec x)^2(\alpha-\beta i)u, \quad \alpha,\beta \in \mathbb{R}
\end{equation}
is used to accelerate the convergence of a Krylov subspace method for solving \eqref{eq:aug}. The preconditioning matrix $M^h$ is obtained from the discretization of \eqref{eq:mu_sl}, in the same way we are discretizing \eqref{eq:helmholtz}. The solution $u$ is computed from the (right) preconditioned system
\begin{equation}
 A^h(M^h)^{-1}\hat{\bfu}^h = \bfg^h, \quad \bfu^h=(M^h)^{-1}\hat{\bfu}^h
\end{equation}
where $A^h$ and $M^h$ are the discrete Helmholtz and shifted Laplacian matrices, respectively. In this paper we focus on the pair $\alpha = 1$ and $\beta= 0.5$, which in \cite{erlangga2006novel} is shown to lead to a good compromise between
approximating \eqref{eq:aug}
and our ability to solve the shifted system using multigrid tools. The SL multigrid method is very consistent and robust for non-uniform models. However, it is considered slow and computationally expensive, especially for large wavenumbers.


\subsubsection*{Geometric multigrid}
Generally, multigrid methods are used to iteratively solve discretized PDEs, like the linear system in \eqref{eq:aug}, defined on a fine grid $\Omega^h$ using a hierarchy of grids. In a nutshell, in a two level setting, we project the fine level problem onto a coarser grid, solve the problem on that coarser grid, and then use the coarse solution to correct the solution on the original fine grid. This is repeated recursively to form a V-cycle, which is applied iteratively to solve the problem.

More explicitly,  multigrid methods are based on two complementary processes: relaxation and coarse grid correction.
The relaxation is obtained by a standard method like Jacobi or Gauss-Seidel, which is only effective at reducing part of the error\footnote{In the case of the Helmholtz system, such relaxation methods do not converge due to the indefiniteness of $A^h$, yet they are effective at smoothing the error using 1-2 iterations only.}. The remaining error, called ``algebraically smooth'', is typically spanned by the eigenvectors of $A^h$ corresponding to small eigenvalues (in magnitude), i.e., vectors $\bfe^h$ s.t.
\begin{equation}\label{eq:aeae}
\norm{A^h\bfe^h}\ll\norm{A^h}\norm{\bfe^h}.
\end{equation}
To reduce such errors, multigrid methods use a coarse grid correction, where the error $\bfe^h$ for some iterate $\bfu^h$ is estimated by solving a coarser system and interpolating its solution:
\[A^H\bfe^H=\bfr^H=I^H_h(\mathbf{g}^h-A^h\mathbf{u}^h),\quad \bfe^h= I^h_H\bfe^H.\]
Here, the matrix $A^H$ is an approximation of $A^h$ on a coarser mesh $\Omega^H$, obtained for width $H = 2h$.
To interpolate the solution from coarse to fine grids we choose the bi-linear interpolation operator $I_H^h$
which is suitable for the problem because the Laplacian operator in \eqref{eq:mu_sl} is homogeneous. The restriction operator $I_h^H$, dubbed ``full-weighting'', projects the fine-level residual to the coarse grid.
Since the coarse problem is still too large to solve, the process of relaxation and coarse grid correction is applied recursively resulting in the V-cycle algorithm given in Alg. \ref{alg:Vcycle}. 

\begin{algorithm}
\SetAlgoLined
 Relax $v_1$ times on $A^h\mathbf{v}^h=\mathbf{f}^h$ with $\mathbf{v}^h$ as an initial guess \;\\
 \eIf{deepest level}{
   $\mathbf{v}^h\leftarrow$ Solve the system $A^h\mathbf{v}^h=\mathbf{f}^h$ directly\;
   }{
   $\mathbf{f}^{2h}\leftarrow I_h^{2h}(\mathbf{f}^h - A^h\mathbf{v}^h)$\;\\
   $\mathbf{v}^{2h}\leftarrow $V-cycle$(\mathbf{v}^{2h}\leftarrow \mathbf{0},\mathbf{f}^{2h})$\;\\
   $\mathbf{v}^h\leftarrow \mathbf{v}^h + I_{2h}^h\mathbf{v}^{2h}$\;
  }
 Relax $v_2$ times on $A^h\mathbf{v}^h=\mathbf{f}^h$  with $\mathbf{v}^h$ as an initial guess\;
 \caption{Multilevel V-cycle: $\mathbf{v}^h\leftarrow$V-cycle$(\mathbf{v}^h,\mathbf{f}^h)$}\label{alg:Vcycle}
\end{algorithm}



In Alg. \eqref{alg:Vcycle} one has to choose the number of levels in the V-cycle. Unlike other scenarios, the smooth error modes of the Helmholtz operator are oscillatory at high wavenumber, and cannot be represented well on very coarse grids. Hence, the performance of the solver deteriorates as we use more levels. For example, the results in \cite{https://doi.org/10.1002/nla.1860} show that the best performance is achieved using three levels only, and the authors suggest to use GMRES as a coarsest grid solver.

\subsection{Machine and deep learning}
The basic purpose of machine learning is to study an approximation of an unknown function $\phi(x)$, using a generic model function $model(x,\theta)$ that has unknown parameters $\theta$. The true function $\phi(x)$ can be, for example, the category of an image $x$, and hopefully the $model$ function is rich enough to approximate $\phi(x)$. The parameters $\theta$, called the weights of the model, are learned
using a set of samples $\{x_i\}$ from a distribution of reasonable inputs (e.g., natural images, spoken words, etc.) such that $model(x_i,\theta)\approx \phi(x_i)$.

In supervised learning, the data set contains pairs $\{ (x_i, y_i) \}_{i=1}^m$ so that $y_i$ is the desired result for the model to provide given $x_i$. To this end, a loss function is defined, and the weights of the model function are estimated by minimizing it. For example, a supervised loss function for predicting the solution of a linear system like \eqref{eq:aug}, can be based on the error mean
\begin{equation}\label{eq:loss}
loss(\theta) = \frac{1}{m}\sum_{i=1}^m \|\mbox{model}(\bfg_i,\theta)-\bfu_i\|_2^2,
\end{equation}
where the training set contains $m$ pairs of $( \bfu_i, \bfg_i)$ satisfying \eqref{eq:aug}.

Deep learning (DL) \cite{GoodBengCour16} is a sub-field of machine learning that is concerned with the collection of model functions belonging to the class of neural networks. A typical neural network architecture includes an input layer, an output layer and a several hidden layers consisting of affine linear transformations and pointwise nonlinear activation functions. The basic layer reads
\begin{equation}\label{eq:layer}
\bfx^{(l+1)} = \sigma(K^{(l)}\bfx^{(l)}+\bfb^{(l)})
\end{equation}
where $\sigma$ is the non-linear activation function, $K^{(l)},\bfb^{(l)}$ are the matrix and bias parameters for the $l$-th layer, respectively. The weights $\theta$ in \eqref{eq:loss} typically include all the parameters for all the layers.
The process of optimizing neural networks---minimizing \eqref{eq:loss}---consists of the forward pass in which the layers like \eqref{eq:layer} are performed one after the other, and the backward pass in which the gradient $\nabla_{\theta} loss$ is computed for applying iterative steps of stochastic gradient descent to minimize the loss. Deep neural networks usually include many hidden layers, that result in a very flexible function representation.

A convolutional neural network (CNN) \cite{krizhevsky2012imagenet} is a type of a deep neural network, which is based on convolution kernels $K$ in \eqref{eq:layer}, leading to a shared-weight architecture to handle data of large images and videos. CNNs are among the most effective methods for dealing with high-dimensional grid-like data, and learning spatial features. They are the leading methods for computer vision tasks like image classification, object detection, and semantic segmentation.

\subsection{Geometric multigrid using CNN modules}\label{sec:MG_GPU}
In our work, both the CNN and the multigrid preconditioner are used interchangeably. We implement the multigrid method using CNN components, in order to enable a good integration of the V-cycle with the CNN in a unified framework within the preconditioner component, using a GPU back-end. The idea is based on the close connection between U-Nets and multigrid V-cycles.

The Helmholtz operator needed for the residual calculation and the Jacobi relaxation is executed using a $3\times 3$ convolution operator and an element-wise vector multiplication with the mass term according to the discretized Helmholtz operator \eqref{eq:ah_matrix} or the discrete \eqref{eq:mu_sl}.
The high order discretization operators in \cite{umetani2009multigrid} can also be obtained similarly.

The geometric multigrid transfer operators, the restriction (coarsening of the mesh) and prolongation (interpolation), are realized by convolution operators with the fixed kernels
\begin{equation}\label{eq:V_Trans_ops}
K_{I_h^H} = \frac{1}{16}\begin{bmatrix}
1 & 2 & 1 \\ 2 & 4 & 2 \\1 & 2 & 1
\end{bmatrix},\quad K_{I_H^h} = \frac{1}{4}\begin{bmatrix}
1 & 2 & 1 \\ 2 & 4 & 2 \\1 & 2 & 1
\end{bmatrix},
\end{equation}
that correspond to the full-weighting and bi-linear operators, respectively. Both kernels are applied with a stride of 2. For the restriction, this results in a coarse mesh $\Omega^H$. For the prolongation,  we use a transposed strided convolution with $K_{I_H^h}$, so that the mesh is refined. The solution for the coarsest mesh problem is done using GMRES with diagonal Jacobi preconditioner \cite{https://doi.org/10.1002/nla.1860}. In addition, all these convolutions are applied with zero padding of size 1, which is consistent with the Dirichlet boundary condition (the ABC is obtained through the mass term). Neumann BC can be obtained by using reflection or replication padding in CNN frameworks,  which result in first and second order Neumann BC, respectively.

\section{Multigrid-augmented deep learning preconditioners}

\begin{figure}
\centering
\includegraphics[width=0.7\textwidth]{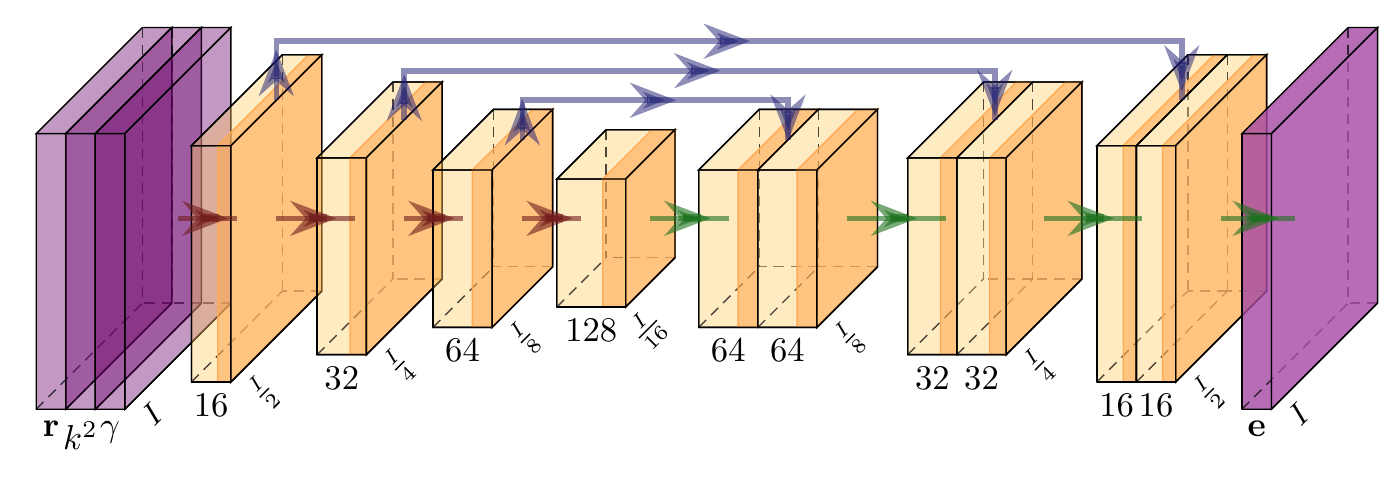}
\caption{Illustration of our U-Net architecture. The purple layers denote the input ($\bfr,\kappa^2,\gamma$) and the output ($\bfe$) of the network, while the orange layers are hidden layers. The red arrows indicate a downsampling convolution operator, followed by a ResNet step, while the green arrows indicate the upsampling convolution operator. The blue arrows indicate feature map concatenations.
}
\label{fig:unet_arch}
\end{figure}

We present a CNN in a U-Net architecture to
accelerate the convergence of the Helmholtz solution, utilizing the GPU capabilities and parallel computation of DL frameworks. In order for the computational complexity to be cost-effective we define the architecture to be relatively small in DL standards (see Fig. \ref{fig:unet_arch} and details later), and yet demonstrate a significant acceleration to the solution. Our DL-based preconditioner is applied within the FGMRES Krylov method \cite{saad1993flexible}. The flexible variant is important as our preconditioner is non-linear, hence also non-stationary.

\subsection{A U-Net architecture} \label{sec:unet_architecture}
The U-Net architecture \cite{ronneberger2015unet} is a CNN for image-to-image tasks such as semantic segmentation, where each pixel is classified. It consists of two rather symmetric parts---contraction and expansion. As described in \cite{ronneberger2015unet}, the purpose of the contraction part is to capture the significant contexts and features, and the purpose of the expansion part is to enables precise localization. In the contraction part, the feature maps are coarsened and the number of channels is increased to get the context vectors.
In the expansion part of the U-Net, the number of feature channels is reduced, but the resolution of the channels is increased. This allows the network to propagate context information from the lower resolution layers to higher resolution layers. Since the contraction and expansion paths are more or less symmetric to each other, this yields a U-shaped architecture.

In our case, given a slowness model and a right-hand-side sampled on a grid (viewed as images), we wish to predict a solution for \eqref{eq:aug} on that grid (an image). Furthermore, when solving the Helmholtz equation, it is important to allow a point source to spread to the entire mesh, as each right-hand-side can be viewed as multiple point-sources. For this reason, a suitable network structure for this task can be a U-Net, since, like a multigrid cycle, it has a multiscale structure. That is, restriction and prolongation layers are performed, allowing each point to affect a wider neighborhood of grid cells (or pixels).   

Our U-Net, presented in Fig. \ref{fig:unet_arch}, receives a complex valued right-hand-side, the slowness model $\kappa^2$, and the attenuation model $\gamma$ (includes the absorbing BC, but can include a heterogeneous attenuation). The non-linearity of the network architecture allows us to learn the complex interactions between these input ``channels'', and minimize the error loss function w.r.t to the weights. In a nutshell, on its way down the hierarchy, the architecture interchangeably applies down-sampling convolutions and ResNet convolution blocks, which are similar in some sense to the Jacobi smoothing steps. Each step applies convolutions and a non-linear activation function. Then, on the way up the hierarchy it applies up-sampling convolutions. The down- and up-sampling are applied using strided convolutions, similarly to  \eqref{eq:V_Trans_ops} with learned weights.

Specifically, on the way down the U-Net, we start with a down-sampling operator, which is performed using a strided convolution (a restriction) layer
\begin{equation}\label{eq:convdown}
  \bfx^{(l+1)} = \sigma (\mathcal{N}(R^{(l)}\bfx^{(l)})) \text{ s.t. } \bfx^{(l)} \in \mathbb{R}^{n_l,n_l,c_l} \text{ and } \bfx^{(l+1)} \in \mathbb{R}^{n_{l+1},n_{l+1},c_{l+1}},
\end{equation}
where $R^{(l)}$ denotes the strided convolution operator, which maps $c_l$ channels to $c_{l+1}=2c_l$ channels while coarsening the channels by a factor of 2 at each dimension ($n_{l+1} = n_{l}/2$). We use a kernel of 5 by 5, so that the network can learn rather high-order transfer operators. $\sigma$ denotes an element-wise exponential linear unit (eLU) 
activation function. We chose a smooth activation function to allow for spatially smooth feature maps which are necessary because our solution is expected to be rather smooth. The batch normalization layers are expressed by $\mathcal{N}$, whose role is to stabilize the optimization of network \cite{batchnorm}.

In between each restriction, while going down the hierarchy, we apply ResNet blocks. The structure of our ResNet block is described as
\begin{equation}\label{eq:resnet}
  \bfx^{(l+1)} = \sigma (\mathcal{N}(\bfx^{(l)} + K_2^{(l)}\sigma(\mathcal{N}(K_1^{(l)} \bfx^{(l)})))).
\end{equation}
The weights $K_1^{(l)}$ and $K_2^{(l)}$ are two convolution operators with $3 \times 3$ sized kernels. ResNet layers are explicitly reformulated as an addition of residual functions with respect to the inputs, instead of learning direct mappings like the standard neural network step in \eqref{eq:layer} or \eqref{eq:convdown}. In \cite{he2015deep}, the authors provide empirical evidence showing that these residual networks are easier to optimize, and can gain accuracy from increased depth. In the ``coarsest'' level of our U-Net, before the prolongation starts, 3 ResNet layers are performed one after the other.

When we go up the hierarchy, we apply the up-sampling convolutions (like prolongations) to refine the channels and decrease the channel space. Such layers are given by
\begin{equation}\label{eq:convup}
  \bfx^{(l+1)} = cat(\mathcal{N}(P^{(l)} \sigma(\bfx^{(l)})), \bfx^{(l')})\\ \text{ s.t. } \bfx^{(l)} \in \mathbb{R}^{n_l,n_l,c_l} \text{ and } \bfx^{(l+1)} \in \mathbb{R}^{2n_l,2n_{l},c_{l}/2},
\end{equation}
where $P^{(l)}$ denotes a strided transposed convolution operator with a 5-by-5 kernel, which maps $c_l$ channels to $c_{l+1}=c_l / 2$ channels. The $cat$ function expresses a concatenation of the current prolonged feature maps $\bfx^{(l)}$ and the corresponding feature maps from the restriction phase $\bfx^{(l')}$. This is marked in arrows in Fig. \ref{fig:unet_arch}, and is similar in spirit to the multigrid's coarse grid correction.

\subsection{U-Net as a preconditioner accelerator}

The U-Net is used to accelerate the solution of the Helmholtz equation. As noted before, since the medium can be very heterogeneous, the solution, even for a point source or a plane wave, can be highly complex at high wavenumber due to reflections and interference. Therefore, we do not pursue a solution via a single application of a network, but aim for a network to be applied several times in the form of a preconditioner. This ensures users that the discrete equation is indeed satisfied to any desired accuracy.

In this work we aim for the U-Net to complement either the Jacobi smoother or the shifted Laplacian multigrid method. Like a relaxation, the SL is efficient at handling error modes with rather large absolute eigenvalues, as the shift does not deviate those eigenvalues by a significant factor. Since SL or a relaxation are efficient at smoothing (=reducing residuals), we aim that the U-Net will focus on the error directly, and it will be accompanied by classical smoothing. That is, we wish to help the U-Net so it can focus on error modes that are not treated well by smoothing or SL. This guides our training procedure, as we describe later on. Furthermore, the smoothing is important since we apply the preconditioner inside a Krylov method, which is sensitive to high residuals. Below we list two preconditioning schemes that we propose.

\subsubsection*{U-Net as a coarse grid solver with Jacobi smoothing}
Our first version of the preconditioner applies one pre- and post- Jacobi relaxations with a U-Net in between. Since the U-Net starts and ends with strided convolutions (like a restriction and prolongation), this ends up being similar to a two level cycle. Alg. \ref{alg:gmres_unet_j} summarizes the process, where the matrix $A^h$ denotes the discrete Helmholtz equation.

\begin{algorithm}[H]
\SetAlgoLined
Algorithm: preconditioner $\bfe\leftarrow M_{\sf JU}(\bfr)$:
\begin{enumerate}
\item Perform a Jacobi step: $\bfe_1 = \mbox{Jacobi}(\bfe_0 = \mathbf{0},\bfr)$.
\item Compute residual, and apply a U-Net:
$
\bfe_2 = \bfe_1 + \mbox{U-Net}(\bfr - A^h\bfe_1).
$
\item Perform a Jacobi step: $\bfe = \mbox{Jacobi}(\bfe_2,\bfr)$.
\end{enumerate}

\caption{ U-Net as a coarse grid solver with Jacobi smoothing.}\label{alg:gmres_unet_j}
\end{algorithm}

\subsubsection*{U-Net as a deflation operator}
A slightly more evolved option would be to apply U-Net, and afterwords to smooth the output using a SL V-cycle. More explicitly:
\begin{equation}
\label{eq:precM_Vu}
\mbox{Preconditioner } M_{\sf VU}(\bfr) = \mbox{V-cycle}(\bfe_0=\mbox{U-Net}(\bfr),\bfr).
\end{equation}
This option somewhat resembles the concept of deflation operator studied in \cite{SHEIKH2016473}, only without the non-linearities, and a coarse grid solution instead of a network. There, the solution is projected using an interpolation operator $Z_{h,H}$, and a coarse grid solution is obtained. That is in addition to SL cycles. Here, the operator corresponding to $Z_{h,H}^T$ is the strided convolution layer at the entrance to the network, and the last U-Net layer is a transposed convolution corresponding to the $Z_{h,H}$ matrix that returns the mesh to its original dimension.



\subsection{Encoder-Solver: an encoder network as a preconditioner setup} \label{sec:encoder}
So far, we introduced a U-Net that receives a residual and a slowness model, and returns an approximate solution. That is, our U-Net generalizes on both the residual and the slowness. This setup is quite challenging compared to the setup in the previous \cite{khoo2018switchnet, moseley2020solving,wiecha2019deep}, where the right hand side and slowness are represented by very few parameters, and the solution is obtained by a single forward pass, and not as a preconditioner. Obviously, the less we have to generalize on, the easier it will be for the network to approximate the solution well. This is also verified in our experiments, where we show that the same network performs better when the slowness model is known during training and during the Krylov solution at test time (see Sec \ref{sec:single_model}).

Having the observation above in mind, we wish to exploit the following: when we solve an equation using several iterations, the slowness model remains the same and is known during those iterations---the U-Net is applied with the same model over and over again, with different residuals. So, in the grand scheme we need to generalize on the slowness  model, but we need not do it at every step. To this end, we propose to split the generalization over the residual and slowness model into two  ``encoder'' and ``solver'' networks. The encoder network receives only the slowness model as input, and prepares feature vectors that include the slowness model information. Those vectors are computed only once at the beginning of the solution process, and are fed into the solver network at every iteration in the solution time. The solver network is applied at every iteration, with the readily available feature vectors inside it. This way we can have a large network for the slowness model, that will be applied only once, and the solver network for the residual can be cheaper. This idea resembles a preconditioner setup that is applied once before the solution process, and is adequate for any right hand side.

\begin{figure}
\centering
\includegraphics[width=0.8\textwidth]{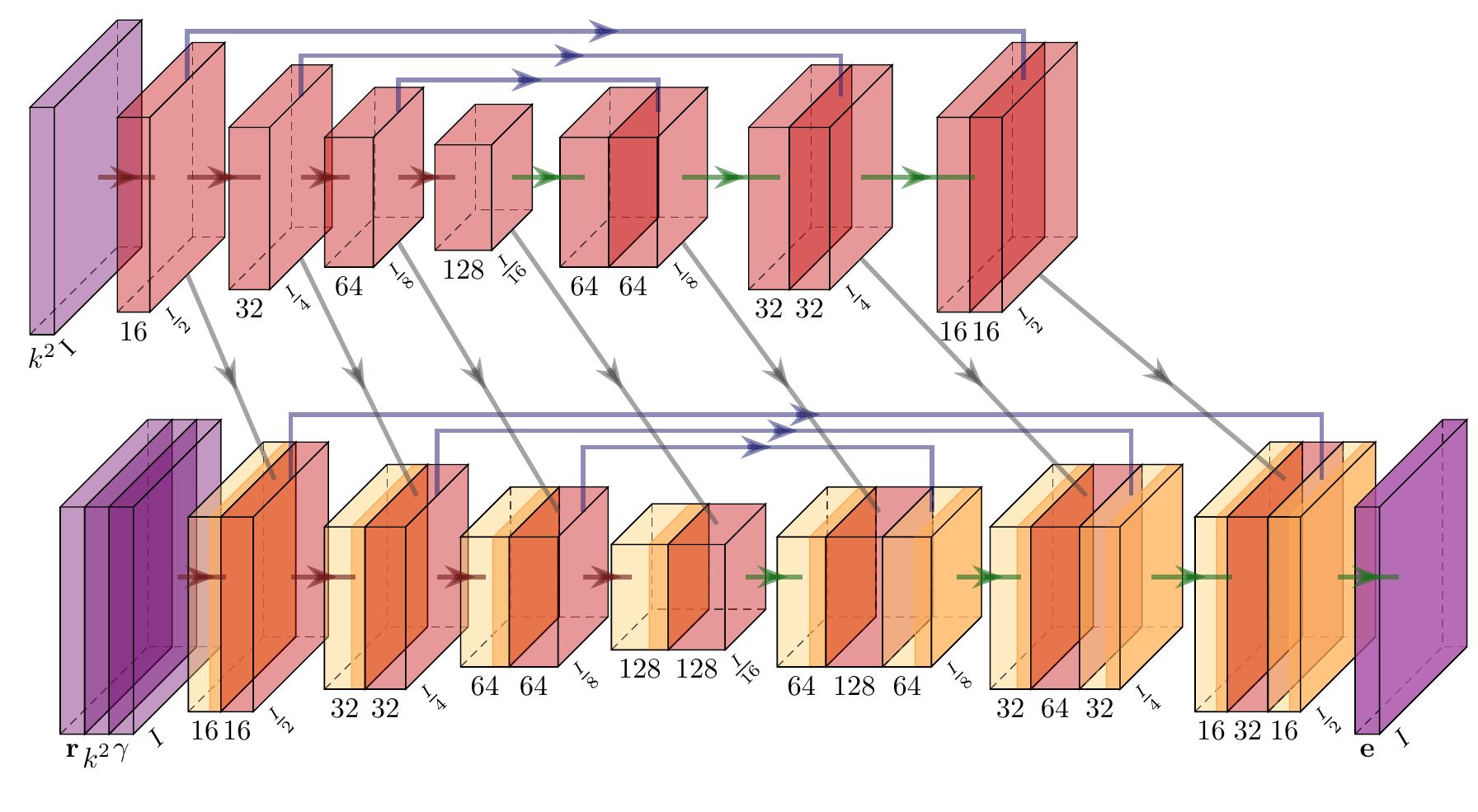}
\caption{Illustration of our Encoder-Solver scheme. The encoder U-Net gets the slowness model and outputs a hierarchy of feature maps. These are fed into the solver U-Net (gray arrows), which receives the residual $\bfr$ and the slowness model $\kappa^2$ as input, and outputs $\bfe$.
}
\label{fig:EncoderSolver}
\end{figure}

Figure \ref{fig:EncoderSolver} illustrates our proposed encoder-solver architecture. The encoder network is also built in a U-Net structure and is almost identical to the network we described in Sec \ref{sec:unet_architecture}. The outputs of the encoder network, the context for the solver network, are hierarchical feature vectors which we concatenate in each layer at the restriction and the prolongation phases of the solver network. The concatenation in the restriction and the prolongation steps can be described as
\begin{equation}\label{eq:concat_hierarchical}
  \bfx^{(l+1)} = \sigma (\mathcal{N}(R_l cat(\bfx^{(l)}, \bff^{(l)})) \text{ and }\bfx^{(l+1)} = cat(\mathcal{N}(P_l \sigma(cat(\bfx^{(l)}, \bff^{(l)}))), \bfx^{(l')})
\end{equation}
respectively, where $\bff^{(l)} \in \mathbb{R}^{n_l,n_l,c_l}$ is the feature vector from the corresponding phase in the encoder network.

\subsection{Training and data augmentation}\label{sec:training_data}
In order for the U-Net to be efficient at solve time as a preconditioner, it must practice or learn on residual vectors that are as close to the solve time FGMRES residuals as possible. Let
\begin{equation}
\bfe^{net} = \mbox{U-Net}(\bfr,\kappa^2;\theta)
\end{equation}
be the forward application of the network for a given right-hand-size $\bfr$ and slowness model $\kappa(\bfx)^2$. As before, $\theta$ denote the collection of network weights (convolutions, biases, etc.) that we wish to learn. Since we wish to complement smoothing processes, in our supervised learning setup we aim to minimize the error with respect to multiple residuals and models. That is, given training data triplets $\{(\bfe_i^{true},\bfr_i,\kappa_i^2)\}_{i=1}^m$, where $A^h\bfe_i^{true}=\bfr_i$, we minimize
\begin{equation}\label{eq:lossnet}
\min_\theta \frac{1}{m}\sum_{i=1}^m\|\mbox{U-Net}(\bfr_i,\kappa_i^2;\theta) - \bfe_i^{true}\|_2^2,
\end{equation}
which is the mean $\ell_2$ error norm, also known as the mean squared error (MSE). The training phase is performed by a stochastic gradient descent optimizer, where in each epoch we sweep through all the training examples in batches and varying order. Specifically, we use the ADAM optimizer \cite{kingma2017adam} with a variable step size.

For supervised learning, we have to provide the network with the corresponding error solution for each residual vector and slowness model. To have such triplets, in some cases one needs to invest a lot of time and resources in preparing the training-set, i.e., solve the PDE many times. However, since we have a linear system here, we can simply generate a random $\bfe_i^{true}$, and easily compute the residual. But, we wish to create a data-set that is as close as possible to the situation in the real solution using FGMRES, and if we simply take the FGMRES residual vectors, we will not have a corresponding error without solving the system. To obtain that we need to generate data with rather smooth vectors in varying ``smoothness levels'', since FGMRES and V-cycles are smoothing operations. Hence, we propose the following procedure for which, in our experience, the training error predicts the final solution's convergence factor quite accurately. First, we generate a random normally distributed $\mathbf{x}_i$ and compute $\bfb_i = A^h\bfx_i$. Then, we apply a random number of FGMRES iterations with SL V-cycle as preconditioner, starting from a zero vector, to get $\tilde\bfx_i$:
\begin{equation}\label{eq:generatexhat}
  \tilde\bfx_i = \text{FGMRES}(A^h, M=\text{V-cycle},\mathbf{b_i},\bfx^{(0)}=\mathbf{0},iter \in \{1,...,10\}).
\end{equation}
Following that, we compute a residual $\bfr_i = \bfb_i - A^h\tilde\bfx_i$ and let $\bfr_i$ be the network input. The true error that we wish the network to recover is given by
\begin{equation}\label{eq:etrue}
\bfe^{true}_i = \bfx_i - \tilde\bfx_i.
\end{equation}
Beyond the procedure above, we insert $\hat{\bfr_i} = h^2 \bfr_i$ as input to the network. This causes the network's input and output to be of the same order of magnitude\footnote{The magnitude of the input and output are important because of how certain default parameters are chosen in DL frameworks---e.g., initialization of weights and biases and scales of activation functions. For this reason, the input is usually normalized, and batch normalization is applied throughout the network.}, and because the system is linear, this change has no effect on the accuracy. Furthermore, to avoid over-fitting and achieve better generalization, we augmented the training examples using a linear combination of the errors and residuals to expand the set. This summarizes the training data generation, given slowness models. The training set and the validation set are manufactured in the same way, and we saw no difference in the loss between these two.

Overall, our network has 4 input channels of the size of the grid. These are the real and imaginary parts of the residuals $\bfr_i$, the real (squared) slowness model $\kappa_i(\bfx)^2$, and the real attenuation model $\gamma(\bfx)$. Although we do not generalize of $\gamma(\bfx)$, we found it slightly beneficial to include it as input, as it allows the network to have a different behaviour at the grid's boundaries.

\begin{remark}
As an alternative to the error loss in \eqref{eq:lossnet}, one may use the residual norm $\|A^h\mathbf{e}_i^{net}-\mathbf{r}_i\|_2$, so that the training can be unsupervised \cite{bar2019unsupervised}. This is a viable option, but it can be a bit misleading, since high errors can in fact have low residuals---these are ``algebraically smooth errors''. Since we can eliminate those errors with V-cycles or Jacobi iterations, we wish to let the network focus on those smooth modes and not on the whole spectrum. We note that empirically, our experience in augmenting the error norms in \eqref{eq:lossnet} with residual norms generally yields inferior performance than focusing on the errors only during training.
\end{remark}


\subsubsection{Slowness models}
Solving the Helmholtz equation for a uniform slowness model is not difficult, as it can be performed using a Fourier transform. When trying to solve this equation over a non-uniform media, the difficulty increases significantly, especially if there are sharp jumps or large ratios between the smallest and highest slowness. Here, we want the network to succeed in different types of models, smooth and non-smooth, that faithfully represent reality. In certain applications (e.g., geophysical or medical imaging) one may have slowness models of a certain typical distribution coming from the application domain. E.g., in geophysical imaging, one typically has a layered model with higher velocity values deeper in the earth. In machine learning tasks, if we wish for an algorithm to generalize on an  input, it is best if we let the algorithm learn on inputs from the same distribution as in the real life application. That is a large part of the advantage of data-driven approaches. Usually, however, to obtain a good generalization we need a large data set of inputs---a collection of slowness models in our case.

\begin{figure}
\centering
    \begin{tabular}{c c c c}
    \includegraphics[width=0.22\textwidth]{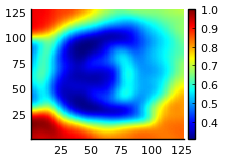} &
    \includegraphics[width=0.22\textwidth]{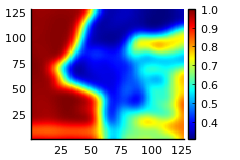} &
    \includegraphics[width=0.22\textwidth]{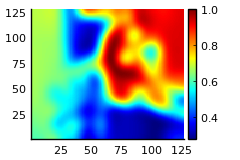} &
    \includegraphics[width=0.22\textwidth]{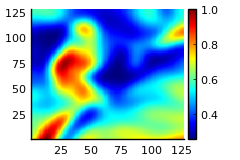}
    \end{tabular}
\caption{Four examples for the $\kappa(\bfx)^2$ models generated from the CIFAR10 image collection. Each model is a random image from the set, enlarged to a grid of 128 by 128 pixels, smoothed using a Gaussian kernel and normalized to the range of [0.25,1].}
    \label{fig:four_models}
\end{figure}

\begin{figure}
\centering
    \begin{tabular}{c c c c}
    \includegraphics[width=0.22\textwidth]{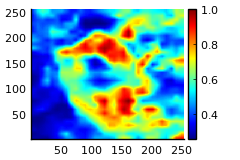} &
    \includegraphics[width=0.22\textwidth]{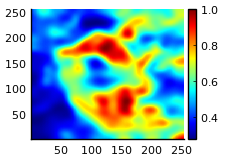} &
    \includegraphics[width=0.22\textwidth]{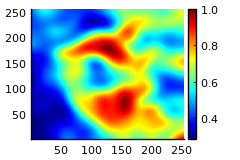} &
    \includegraphics[width=0.22\textwidth]{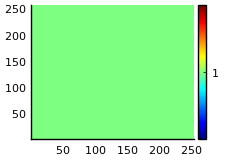} \\

    \includegraphics[width=0.22\textwidth]{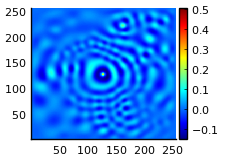} &
    \includegraphics[width=0.22\textwidth]{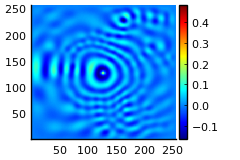} &
    \includegraphics[width=0.22\textwidth]{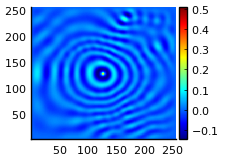} &
    \includegraphics[width=0.22\textwidth]{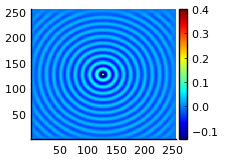}    \\

    (a) & (b) & (c) & (d)
    \end{tabular}
    \caption{Illustration of smooth levels of $\kappa(\bfx)^2$, and their effect on the full solution of the Helmholtz equation for a point source. In the upper row we show the model $\kappa^2$ for each test. The leftmost, (a), is an image from the CIFAR10 data-set, enlarged to a $256\times256$ grid, and the two inner figures, (b) and (c), are the same model after a Gaussian smoothing with kernel sizes of 5 and 10, respectively. The rightmost, (d), is the constant homogeneous model. On the bottom row we show the full solution for the corresponding model in the first row.}
    \label{fig:smooth_models}
\end{figure}
\begin{figure} \centering
    \begin{tabular}{c c c c}
    \includegraphics[width=0.22\textwidth]{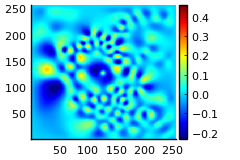} &
    \includegraphics[width=0.22\textwidth]{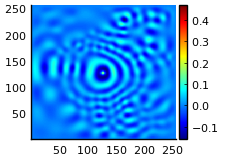} &
    \includegraphics[width=0.22\textwidth]{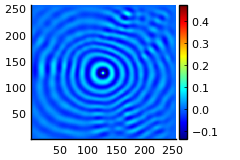} &
    \includegraphics[width=0.22\textwidth]{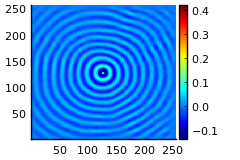} \\
    (a) & (b) & (c) & (d)
    \end{tabular}

    \caption{Illustration of the effect of the value ranges of $\kappa^2$ on the full solution of the Helmholtz equation for the point source problem. The leftmost, (a) is a full solution for kappa normalization between 0 and 1, (b) $\kappa^2 \in [0.25,1]$, (c) $\kappa^2 \in [0.5,1]$ and (d) $\kappa^2 \in [0.75,1]$.}
    \label{fig:range_models}
\end{figure}

Here, we use a rather general slowness model collection to demonstrate our method. The models that we use for the U-Net training, are created from a large set of natural images CIFAR-10 \cite{Krizhevsky09learningmultiple}. This data set is usually used as image classification benchmark and contain 50K images of resolution $32\times 32$ pixels. We use this data set and manipulate the images to become slowness models in order to have a significant amount of training models to learn from. In particular, we use a training set of $20,000$ random sample images.
To transform the natural images into slowness models, we first enlarge them to the grid size using a bi-linear interpolation, then smooth them with a single application of a Gaussian smoothing kernel, and normalize the values to the range of $[0.25,1]$. Fig. \ref{fig:four_models} shows four different models produced from the CIFAR-10 image collection.
Fig. \ref{fig:smooth_models} presents how the smoothness level of the model, $\kappa^2$, affects the difficulty of the solution. It is evident that as the medium is less smooth, the solution is more complex and is harder to model by a neural network. Similarly, the range of slowness values in the model also affects the the complexity of the solution. Fig. \ref{fig:range_models} illustrates how the solution is again more complex as the contrast in the slowness model is larger.

\subsection{CNNs vs. PINNs and a proposed mini re-training}

In the approach of PINNs, one approximates the true solution $u(x)$ as a point-wise neural network $u_{nn}(x)$ that receives $x\in\mathbb{R}^d$ in a mesh-free manner, and possibly another simple property such as a source location, or a plane-wave direction for the right-hand-side. The PDE is then solved by the optimization process for learning the neural network's weights, and any underlying heterogeneous medium or boundary conditions are assumed to be incorporated into the learned weights through the loss minimization.
However, this way, we are quite limited in generalizing over the model coefficients, and must have a powerful enough architecture for $u_{nn}$ that is able to approximate the solution $u(x)$ well enough with some learned weights. The upside is that the optimization hopefully leads to the best approximation possible for the specific solution given the architecture. As the learning is in fact the solution of the PDE, we also need to assume that the optimization is not expensive and does not require the forward-backward passes over too many points.

In some sense, if we define a CNN with $1\times1$ convolutions only, and feed it with the function $x$ on a grid or a mesh, then we get multiple separate instances of $u_{nn}$ in each forward pass for all the points on the grid. A real CNN will also incorporate spatial connections between neighboring (grid) points, hence it is a more powerful model. Similarly, in cases where a spatially adaptive mesh is more suitable than a regular grid, a graph neural network can replace the role of the CNN for the same tasks \cite{eliasof2020diffgcn}. So - if we would define our CNN with $\bfx$ as input, and consider training in solve time, we get a CNN version of PINNs. The true difference, however, is that here, we try to generalize over grid-sized values, like $\bfr$ and $\kappa^2$, which is the way CNNs are used in computer vision, for example.

With the comparison above in mind, our next proposed approach aims to further improve the solver U-Net for the slowness model by a mini retraining procedure at solve time for the given slowness model. That is, given a slowness model $\kappa^2$ we apply a few training iterations over a randomly generated data set of residuals to improve the U-Net weights. This is similar to the process obtained in PINNs, and is motivated by our relatively better performance for known models compared to the case where we generalize over the models. This option would fit cases where we need to solve multiple linear systems with the same slowness model but for many right-hand-sides. Such a scenario fits the case of inverse problems, where the medium needs to be estimated. There, the slowness model evolves iteratively and does not change a lot between consecutive iterations of the inverse solution. Hence, the model retraining can gradually improve the forward solver network with the iterations, like in the PINNs forward/inverse solver in \cite{bar2019unsupervised}. We note that here our retraining phase uses the error-based supervised training procedure described earlier, and not the residual-based training in the common approaches of PINNs.

\section{Numerical Results}
In this section we evaluate the performance of the preconditioner together with FGMRES for solving the Helmholtz equation for random right-hand-sides. In particular, we use a block restarted FGMRES(10) described in \cite{calandra2012flexible}, with a subspace of size 10, solving 10 different right-hand-sides together (a block size of 10). We consider it to be the preconditioner test, and later discuss and demonstrate the influence of the number of right-hand-sides  and Krylov subspace size. In addition, we present the relative error norm plot throughout the training period, which shows the value of the loss function \eqref{eq:lossnet} vs. the epoch number. Ideally, the asymptotic value of the relative training error should indicate the real convergence factor at solve time using FGMRES. It is indeed approximately so thanks to the varying smoothing levels that we apply for our training residuals in Eq. \eqref{eq:generatexhat}. Using this plot, we can also estimate over-fitting and adjust the hyper-parameters (learning rate and mini-batch size, for examples). Our code is written in the Julia language \cite{Julia-2017}, utilizing the Flux DL framework \cite{Flux2018}.

Unless stated otherwise, the test cases that we present are for \eqref{eq:helmholtz} defined on a 2D domain of $[0,1]^2$, discretized on a nodal regular grid of $128 \times 128$ cells, and a wave frequency defined by $\omega=20\pi$, which leads to about 12.5 grid points per wavelength. We use an absorbing boundary layer of width 10 to prevent reflection from all the model boundaries. We compare three preconditioners: The SL V-cycle $M_{\sf V}$, $M_{\sf JU}$ in Alg. \ref{alg:gmres_unet_j} and $M_{\sf VU}$ for Eq. \eqref{eq:precM_Vu}. We evaluate the methods as preconditioners by their convergence factor, calculated as
\begin{equation}\label{eq:convergenceFact}
\rho_{\sf M} = \sqrt[T]{\|\bfr^{(T)}\|_2 / \|\bfr^{(0)}\|_2},
\end{equation}
where $T$ is the number of iterations required to reduce the residual norm by a factor of $10^6$. For the SL V-cycles we use a shift of $\beta=0.5$, 3 levels, 1 pre- and post- damped Jacobi smoothing with a damping factor of 0.8, and a coarsest grid solution of GMRES(10) with the Jacobi diagonal preconditioner as suggested in  \cite{https://doi.org/10.1002/nla.1860}. Under these settings, the convergence factor of FGMRES using the SL V-cycle preconditioner for a model $\kappa^2 \in [0.25,1]$ is
$\rho_{\sf V} \approx 0.9$,
and more than $130$ V-cycles are required to reach convergence. The average convergence factor of $M_{\sf V}$ for a grid of $256 \times 256$ cells, and frequency $\omega=40\pi$, is $0.974$. The advantage of SL is that $M_{\sf V}$ is quite consistent and robust for uniform, smooth and non-smooth models---the difference in the convergence factors is $\pm 0.03$ only, depending on the model. Note that the convergence factor according to \eqref{eq:convergenceFact} is typically lower than the asymptotic convergence factor.

The training process for all the U-Nets is performed by the ADAM optimizer \cite{kingma2017adam}, and includes 80-120 epochs. The initial learning rate is $10^{-4}$ and the initial mini-batch size is 20. 
Along the epochs we reduce the learning rate and increase the mini-batch size. The schedule of these changes as well as the amount of epochs we apply vary a bit between the different experiments and are adjusted by cross validation on the test set. In the majority of the cases, after $t_0=48$ epochs we divide the initial learning rate by 10 and double the mini-batch size. The next changes of the learning rate and mini-batch size occur after every $t_i = \frac{t_0}{2^i}$ epochs. That is, the changes occur at epocs 48, 72, 84, 90 etc, until the learning rate essentially reaches zero. These settings can be used as default parameters for all the experiments that we show.



\subsection{Performance of a homogeneous model}
\begin{figure}
\centering
    \begin{tabular}{c c}
        \begin{tabular}{c c}
            \scriptsize{(a)} & \scriptsize{(b)}\\
            \includegraphics[scale=0.38]{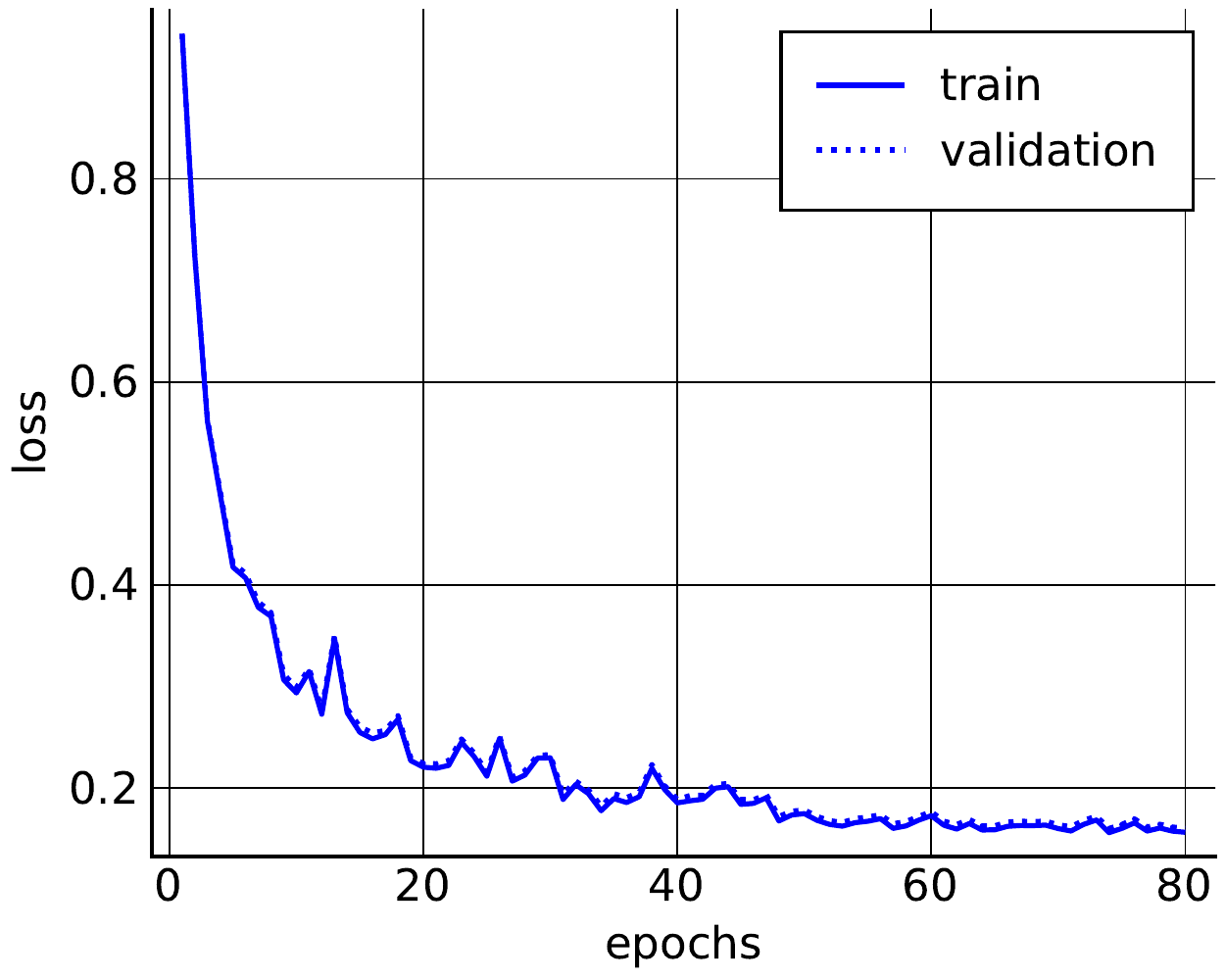}  &
            \includegraphics[scale=0.38]{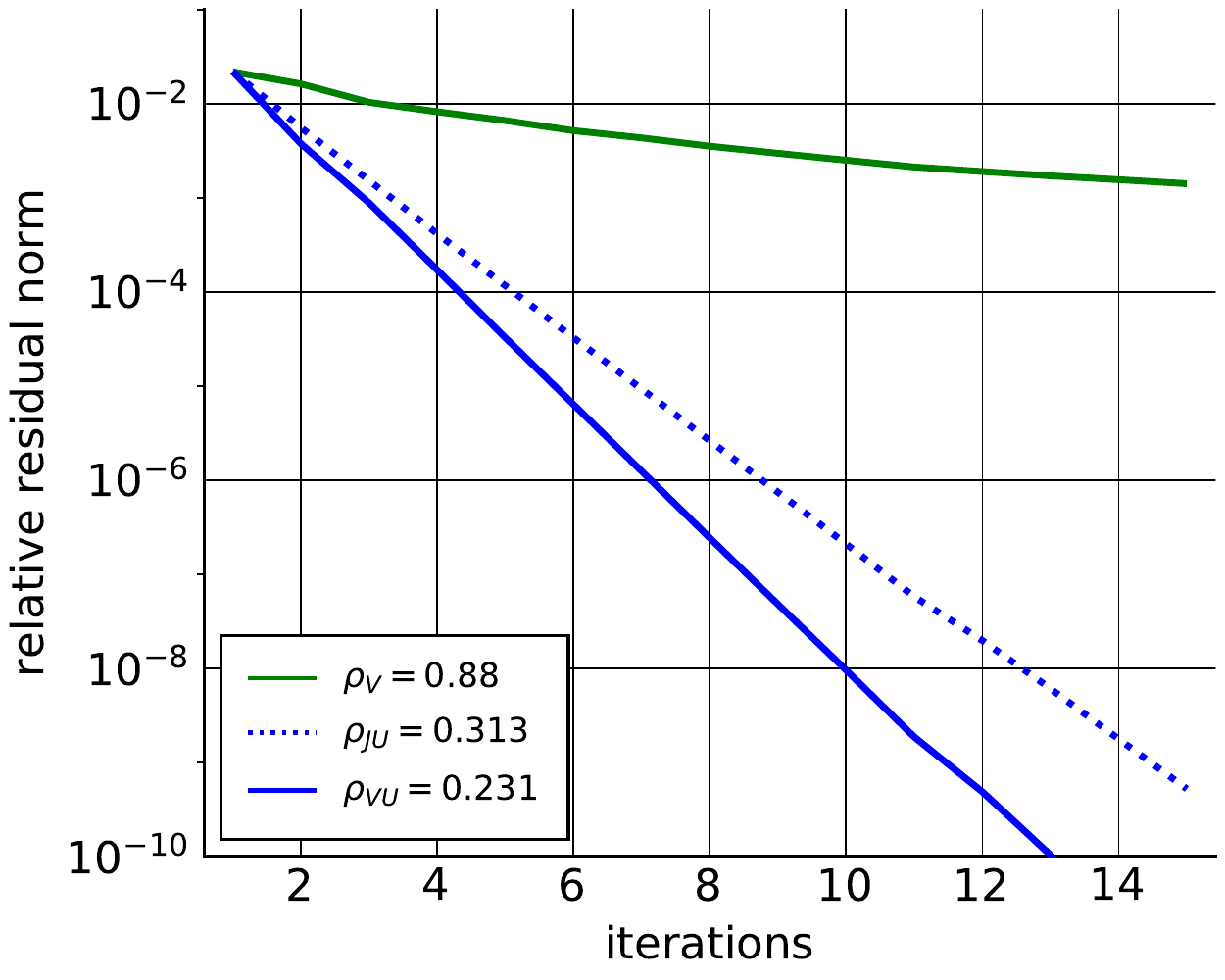}
        \end{tabular} &
        \begin{tabular}{c}
            \scriptsize{(c)} \\
            \includegraphics[scale=0.46]{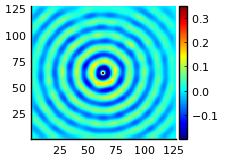} \\
            \scriptsize{(d)}\\
            \includegraphics[scale=0.46]{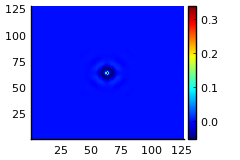}
        \end{tabular}
    \end{tabular}
    \caption{Experiments for a homogeneous slowness model. (a) the training plot. (b) the residual drop history of the FGMRES method.  $\rho_V$,$\rho_{\sf JU}$,$\rho_{\sf VU}$ are the convergence factors of the preconditioners $M_{\sf V}$, $M_{\sf JU}$ and $M_{\sf VU}$,  respectively. (c) shows the real part of the result of a single U-Net operation on a point source, compared to the result of a single V-cycle for the same input in (d).}
    \label{fig:homogeneous_results}
\end{figure}

In the homogeneous test case, we solve the system \eqref{eq:aug} for a uniform model, known at training. That is, the training is performed on random residual and error vectors only as described earlier. Fig. \ref{fig:homogeneous_results} describes the results of the proposed method for this case. The average relative error, i.e. the minimum estimate value of the loss function \eqref{eq:lossnet}, of the training and test sets at the end of the training process is $0.18$. The FGMRES method with U-Net-based preconditioner as described in Eq. \eqref{eq:precM_Vu} converges within 10 iterations with convergence factor $\rho_{\sf VU} = 0.231$ compared to $0.88$ using only V-cycles. Our cheaper preconditioner, which uses only our U-Net and Jacobi steps as described in Alg. \ref{alg:gmres_unet_j}, converges within 12 iterations with convergence factor $\rho_{\sf JU} = 0.313$. These are very favorable convergence factors given how poorly V-cycles behave for a homogeneous model. On the other hand, homogeneous problems can be efficiently solved using fast Fourier transform.

\subsection{Performance on a single heterogeneous model} \label{sec:single_model}
\begin{figure}
\centering
    \begin{tabular}{c c c c}
    \multicolumn{2}{c}{\scriptsize{(a)}} & \multicolumn{2}{c}{\scriptsize{(b)}}\\
    \multicolumn{2}{c}{\includegraphics[scale=0.38]{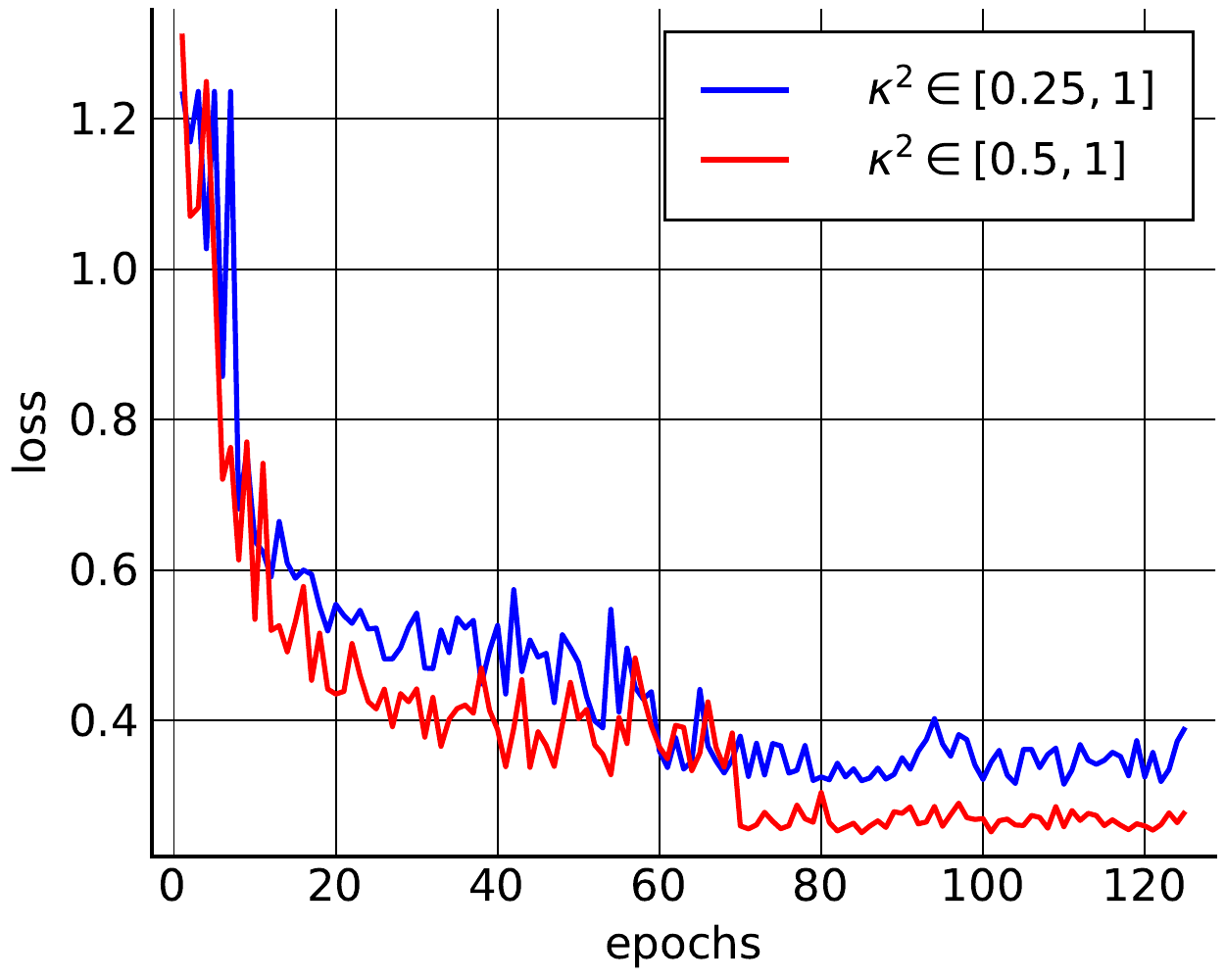}} & \multicolumn{2}{c}{\includegraphics[scale=0.38]{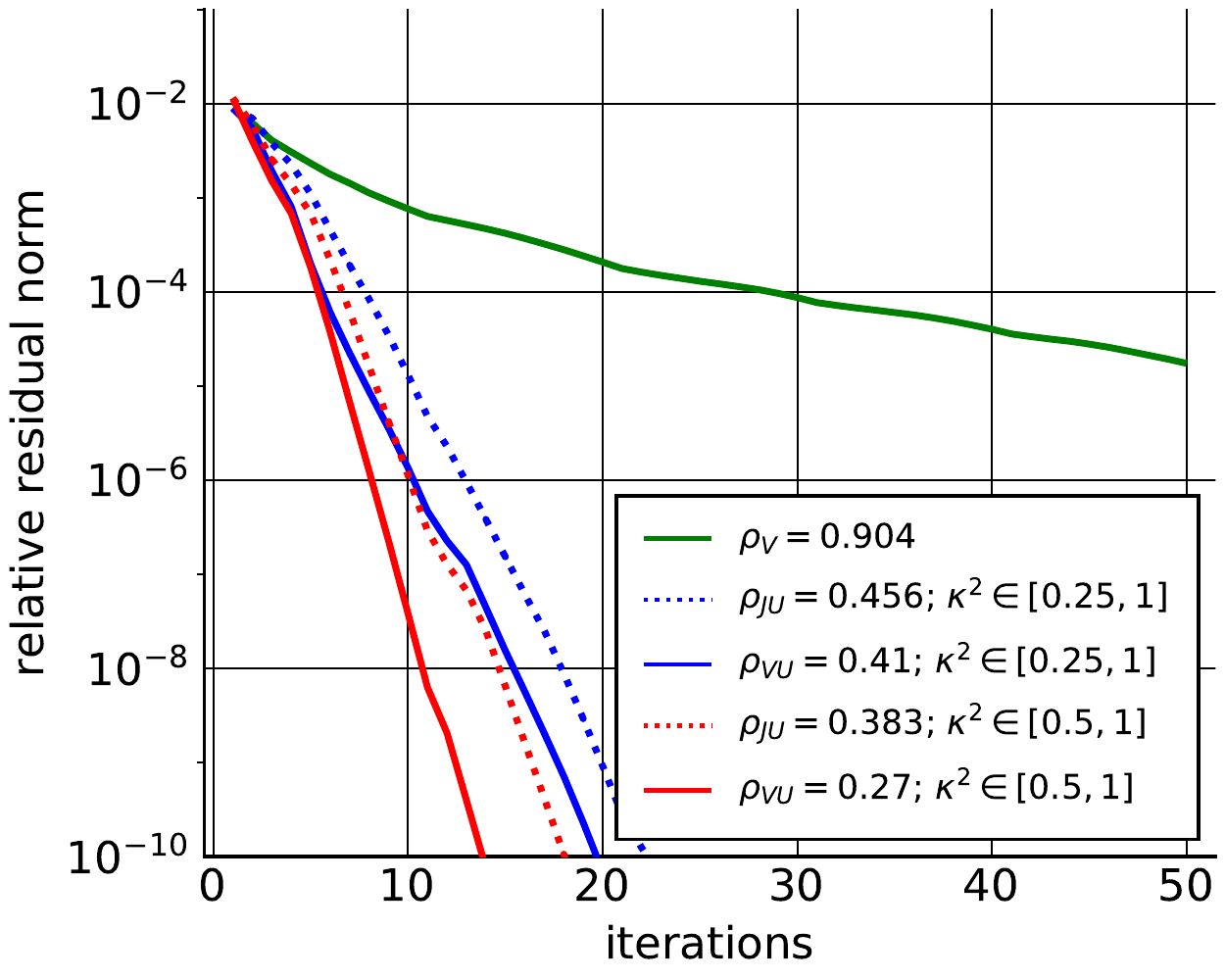}}\\
    \scriptsize{(c)} & \scriptsize{(d)} & \scriptsize{(e)} & \scriptsize{(f)}\\
    \includegraphics[scale=0.48]{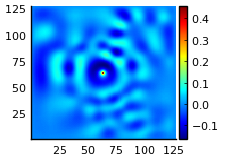} &
    \includegraphics[scale=0.48]{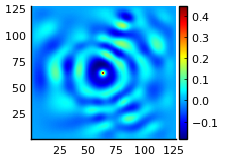} &
    \includegraphics[scale=0.48]{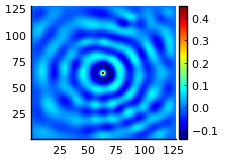} &
    \includegraphics[scale=0.48]{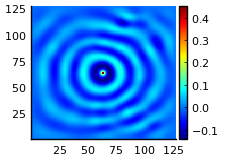}
    \end{tabular}
    \caption{Experiments for a known heterogeneous slowness model. We compare models for which $\kappa^2 \in [0.5,1]$, and $\kappa^2 \in [0.25,1]$. (a) shows the convergence in the training phase and (b) shows the results of the three preconditioner tests. The bottom figures describe the real part of the results of a single U-Net operation, (c) for $\kappa^2 \in [0.25,1]$ compared to the real part of the full solution in (d), and (e) for the case of $\kappa^2 \in [0.5,1]$ compared to the full solution (f).}
    \label{fig:single_results}
\end{figure}

In this section we examine our method for a single heterogeneous slowness model, known during the training phase. In this scenario, the network can learn weights specific to that model. As shown in Fig. \ref{fig:single_results}, for a model normalization of $\kappa^2\in[0.5, 1]$ we obtained convergence factors of $\rho_{\sf VU} = 0.27$ and $\rho_{\sf JU} = 0.383$, which are close to the performance on the uniform model.
For the more difficult case of $\kappa^2 \in [0.25,1]$, we obtain $\rho_{\sf VU} = 0.41$, $\rho_{\sf JU} = 0.456$. Furthermore, in plot (c) of Fig. \ref{fig:homogeneous_results} and in the bottom plots (c) and (e) of Fig. \ref{fig:single_results}, we show the results of a single iteration of Eq. \eqref{eq:precM_Vu} on a centralized point source, a problem that is not included in the training set. It is evident that U-Net has indeed learned to produce a global solution, and the integration with the V-cycle is intended to ``clean'' or smooth the solution and ensure convergence.

\subsection{Generalization on the slowness model}
As we showed in the previous section, our U-Net can learn to solve \eqref{eq:aug} quite efficiently for a given heterogeneous slowness model. However, the training of the U-Net is quite expensive to apply at solve time. In a more applicable approach for the realistic world, we wish the neural network to work and generalize on models it had not seen in the training phase. Therefore, we will now present the ability of the proposed method to generalize for different models. In addition, we will also present the upgraded versions of our approach: the encoder-solver dual networks and the mini re-training approach.

In the experiments in this section, the training phase was performed on a set of $20,000$ pairs of $( \bfr,\bfe)$ produced as described earlier in Sec. \ref{sec:training_data}, and $\kappa^2$ were randomly generated from a collection of the $50,000$ images in the CIFAR10 data-set. In the preconditioner test, we compare the performance of our methods to that of V(1,1)-cycle, which is also used for $M_{\sf UV}$. We will examine the capabilities of the preconditioner using
\begin{enumerate}
    \item A stand-alone U-Net with smoothing/deflation.
    \item The encoder-solver framework.
    \item All options of the U-Nets after a modest ``mini retraining'' phase at solve time.
\end{enumerate}

\begin{figure}
\centering
    \begin{tabular}{c c c}
        \scriptsize{(a)} & \scriptsize{(b)} & \scriptsize{(c)}\\
        \includegraphics[scale=0.35]{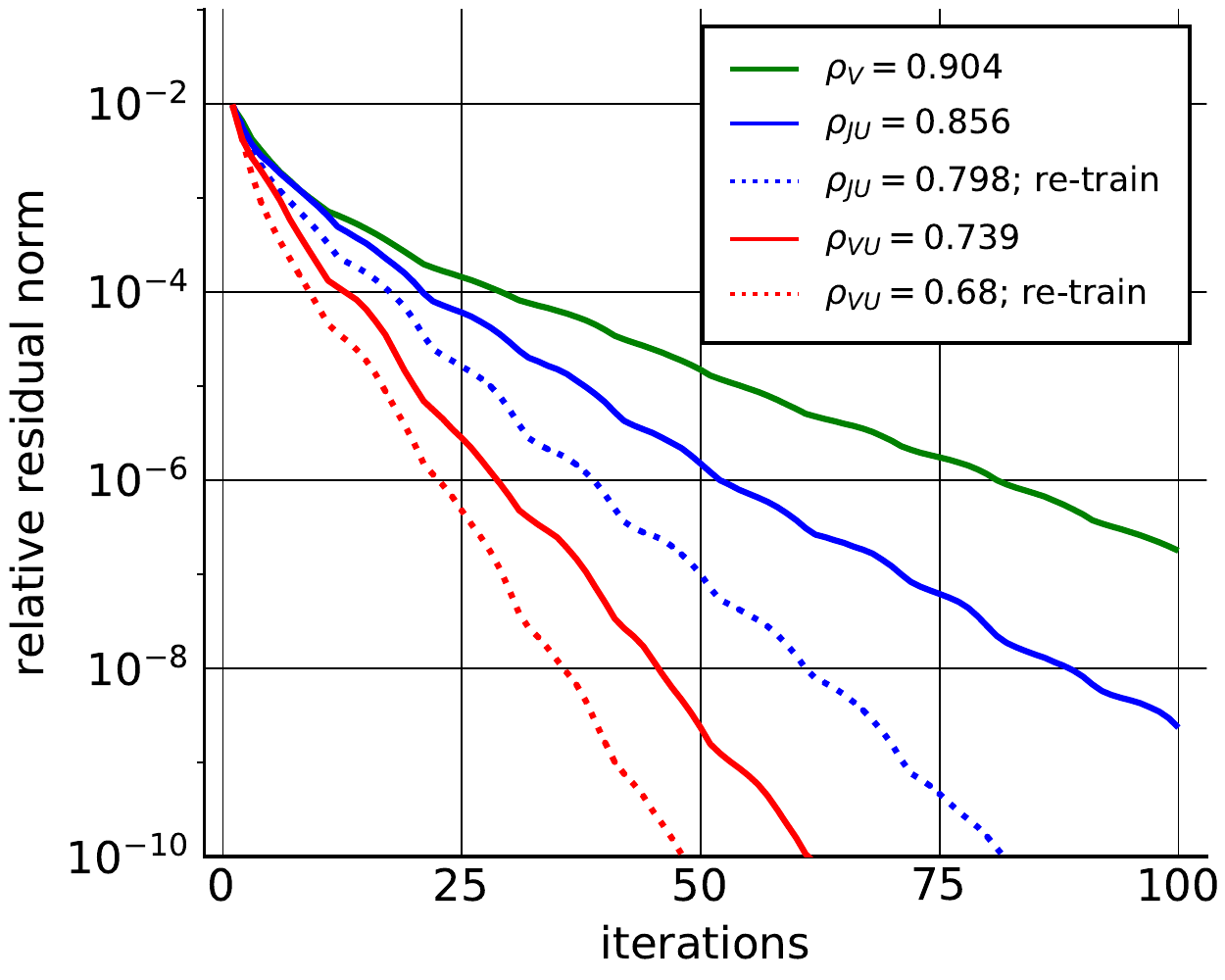} & \includegraphics[scale=0.35]{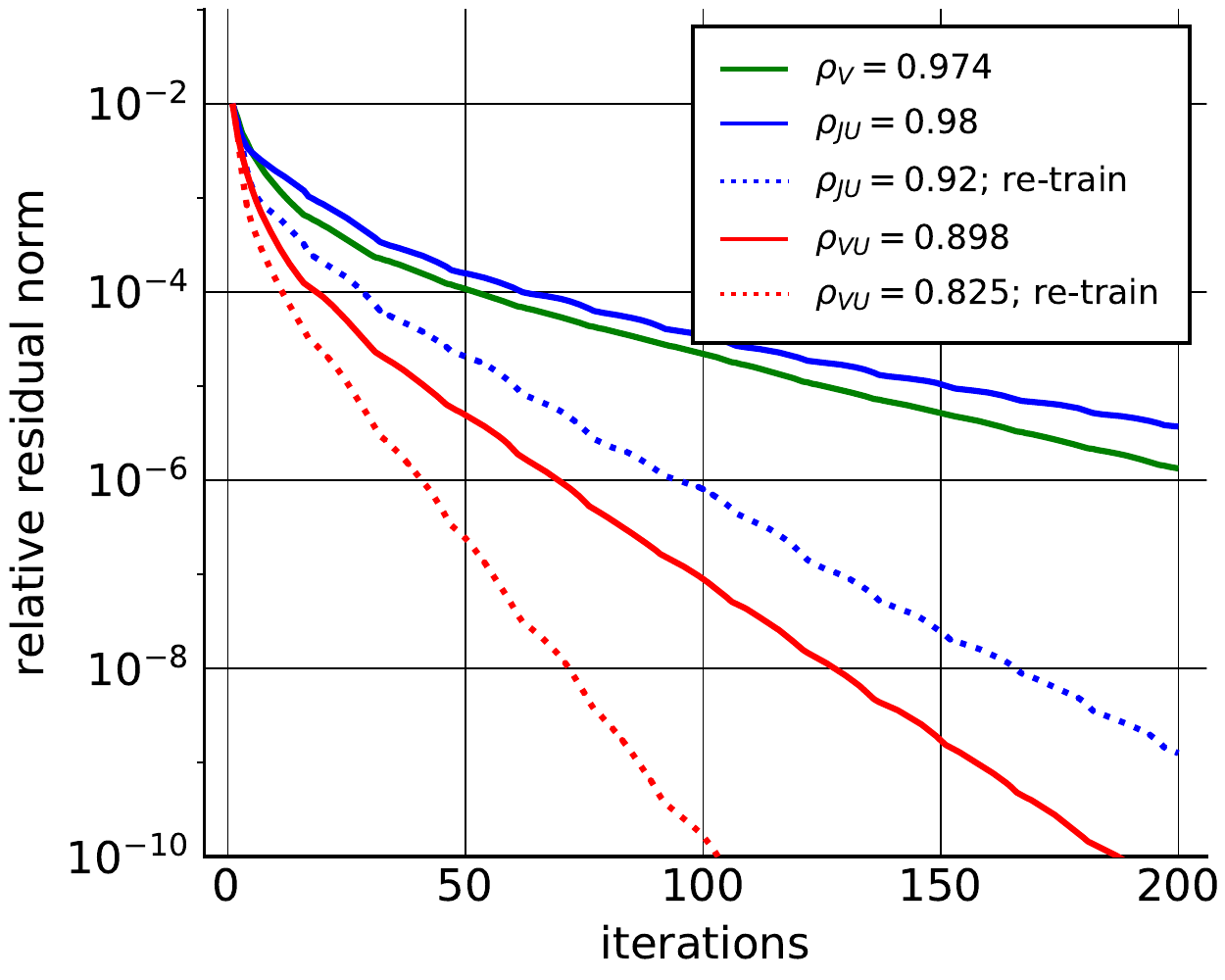} &
        \includegraphics[scale=0.35]{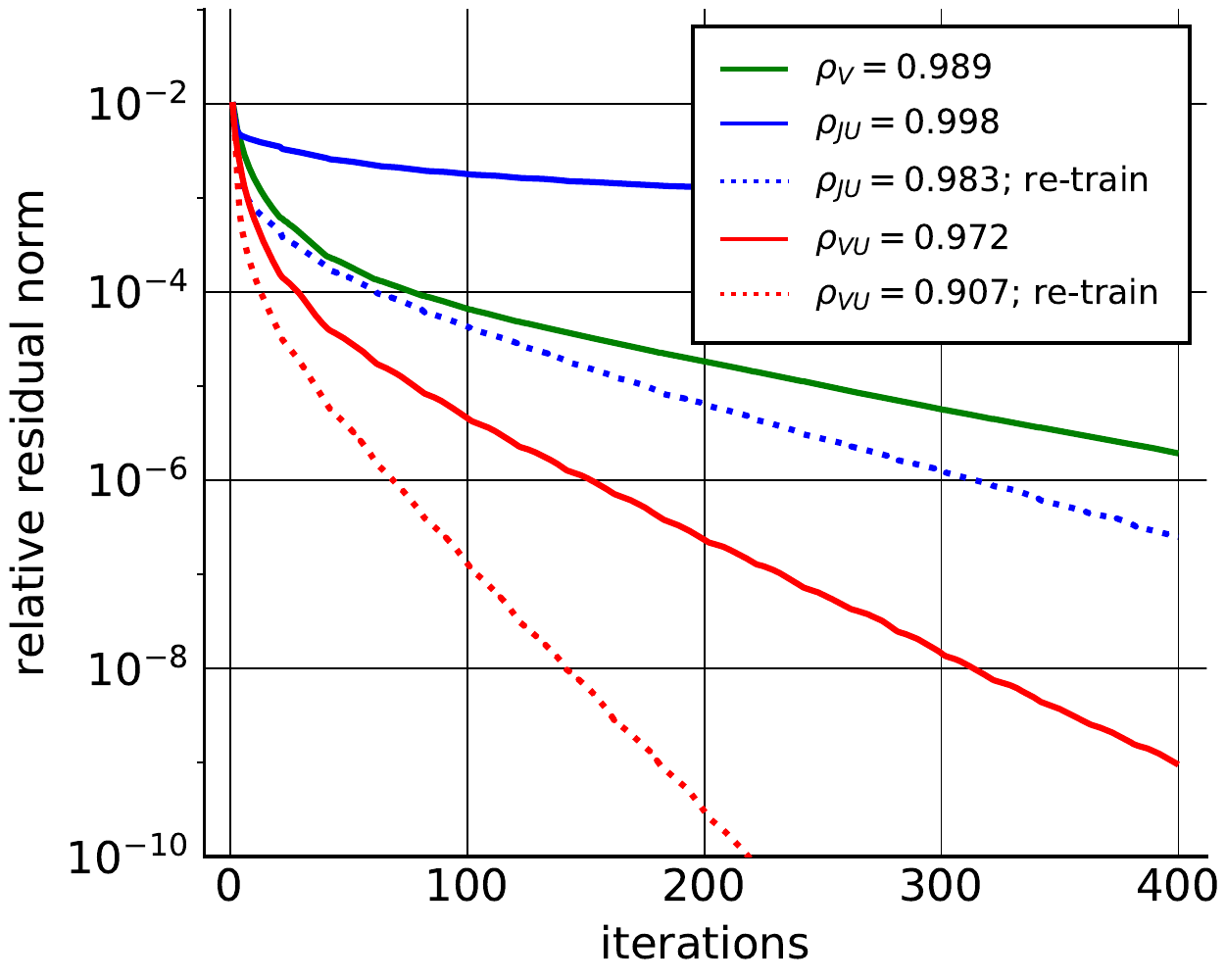} \\

         \scriptsize{(d)} & \scriptsize{(e)} & \scriptsize{(f)} \\
          \includegraphics[scale=0.35]{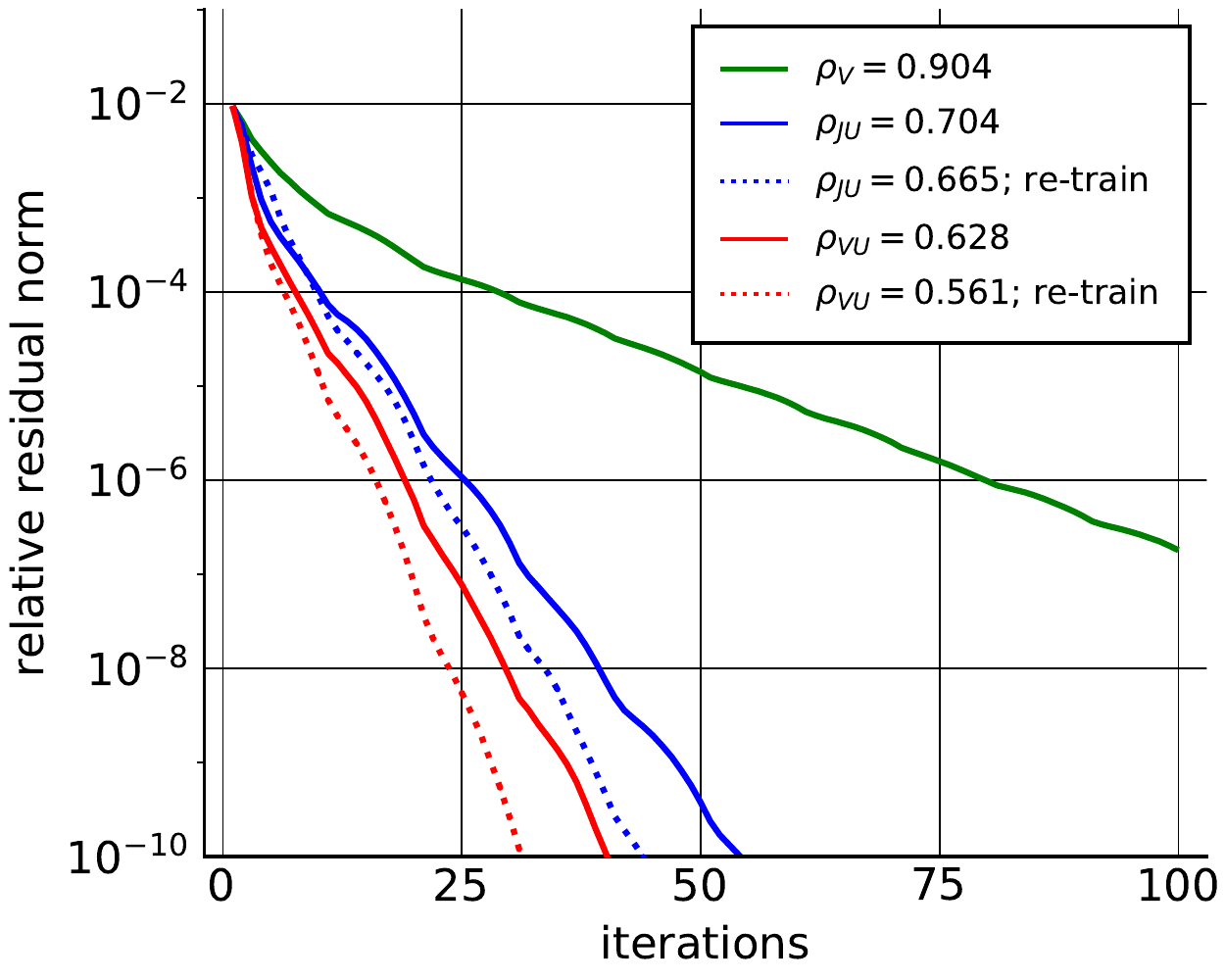}  & \includegraphics[scale=0.35]{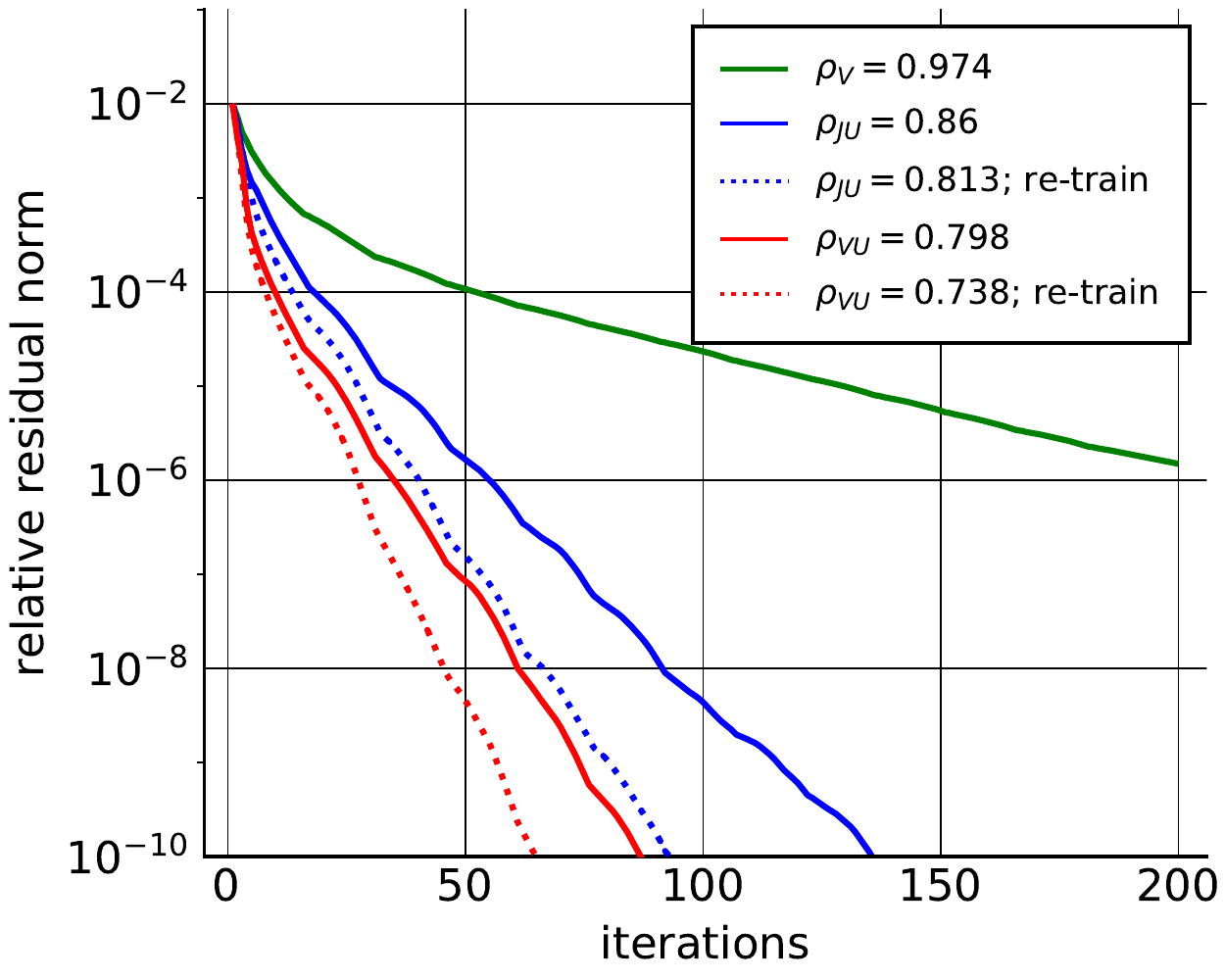} &
          \includegraphics[scale=0.35]{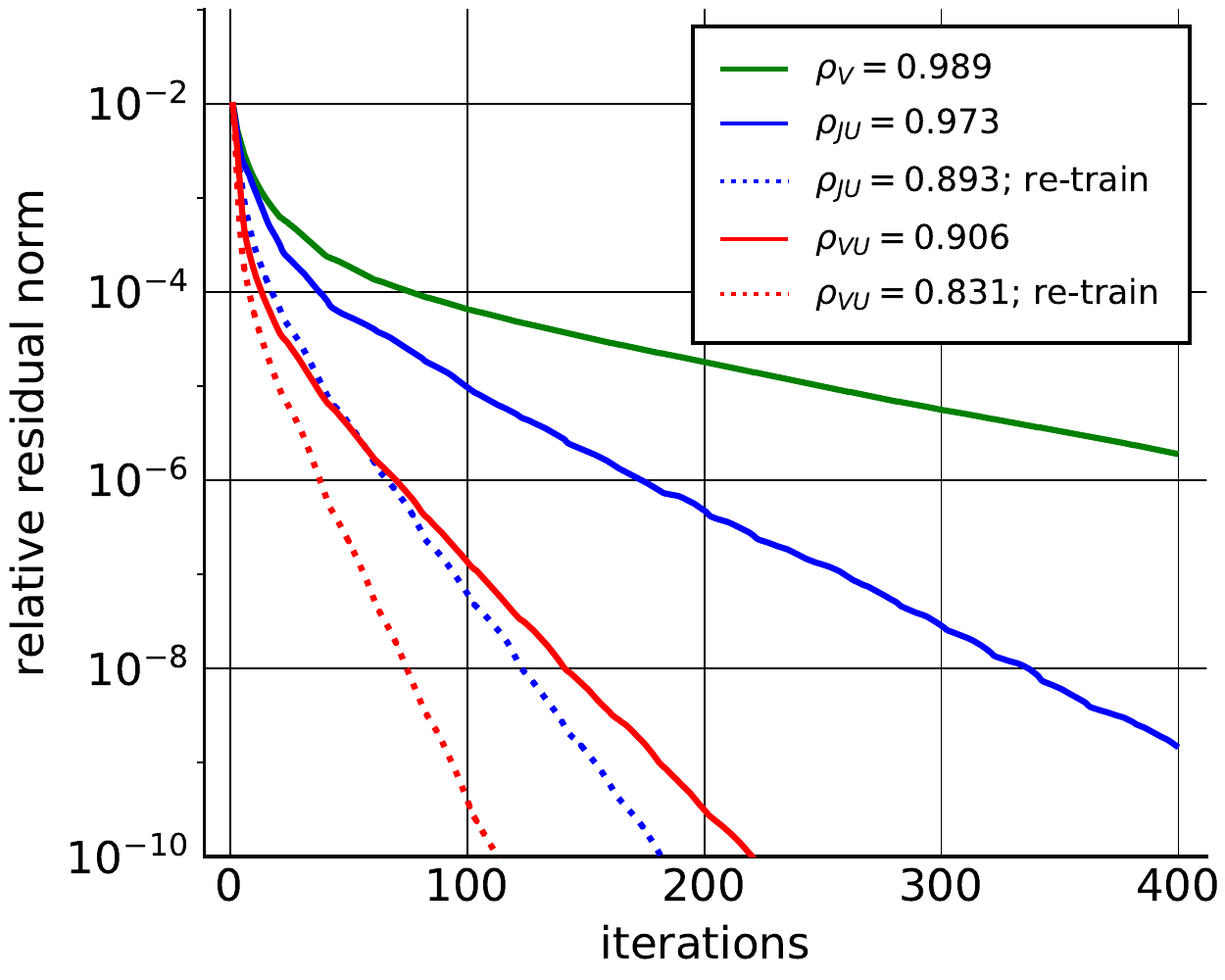}
    \end{tabular}
    \caption{The convergence history of the solutions of \ref{eq:aug} for random RHSs.
    First row: stand-alone U-Net for problem of size (a) $128\times 128$ , (b) $256 \times 256$ and (c) $512 \times 512$. Second row: the encoder-solver U-Nets the same experiments.
    Solid and dotted lines mark the  preconditioners without and with the mini-retraining phase at solve time, respectively. 
    }
    \label{fig:retrain_results}
\end{figure}

Table \ref{table:convergence_factor_table} and Figure \ref{fig:retrain_results} summarize and compare the results for all three cases, and in the next subsections we discuss each case separately. For all options we examine how the solver performs in the case of its training setting. That is, we train the network on a grid of $128\times128$ for $\kappa^2\in[0.25,1]$, but also examine the performance on larger grids and frequencies, and for $\kappa^2\in[0.5,1]$. In particular, we solve \eqref{eq:aug} for random right-hand-sides on grids of size $256 \times 256$ and $512 \times 512$ with $\omega = 40\pi,80\pi$, respectively , using a network trained on residuals and slowness models of size of $128 \times 128$ only.
Table \ref{table:convergence_factor_table} presents the results for the encoder-solver setting. Regardless of the specific preconditioner, we see that if the networks are trained on a single model, they exhibit the best performance. However, they fail to generalize on slightly different cases---even for the easier normalization of the same model (row (a) vs row (b) in the single model column). If we train the networks to generalize over slowness models, it can also generalize on the problem size and slowness contrast. All cases are improved using mini-retraining.

\begin{table}[]
    \centering
    \begin{tabular}{|c|c|c|cc|cc|cc|}
        \hline
     \multicolumn{2}{|c|}{} & & \multicolumn{2}{c|}{Single model  }
        & \multicolumn{2}{c|}{General model  } & \multicolumn{2}{c|}{Re-trained   } \\

        \multicolumn{2}{|c|}{} & $M_{\sf V}$ & $M_{\sf VU}$ & $M_{\sf JU}$& $M_{\sf VU}$ & $M_{\sf JU}$& $M_{\sf VU}$ & $M_{\sf JU}$ \\

         \hline
         \multirow{2}{4cm}{(a) $\kappa^2 \in [0.25,1];n=128$}
         & $\rho$ & $0.904$ & $0.41$ & $0.456$ & $0.628$ & $0.704 $ & $0.561$ & $0.665$ \\
          & $T$ & $137$ & $16$ & $19$ & $30$ & $40$ & $24$ & $34$ \\
         \hline
         \multicolumn{9}{|c|}{Generalization for cases not seen in training}\\

         \hline
         \multirow{2}{4cm}{(b) $\kappa^2 \in [0.5,1];n=128$}
         & $\rho$ & $0.867$ & $0.827$ & \multirow{2}{0.8cm}{fail} & $0.516$ & $0.573$ & $0.432$ & $0.498$ \\
        & $T$ & $97$ & $73$ &  & $21$ & $25$ & $17$ & $20$ \\
         \hline
         \multirow{2}{4cm}{(c) $\kappa^2 \in [0.25,1];n=256$}
         & $\rho$ & $0.974$ & $0.958$ & \multirow{2}{0.8cm}{fail} & $0.798$ & $0.86$ & $0.743$ & $0.811$ \\
         & $T$ & $391$ & $240$ &  & $62$ & $92$ & $46$ & $67$ \\
         \hline
         \multirow{2}{4cm}{(d) $\kappa^2 \in [0.5,1];n=256$}
          & $\rho$ & $0.946$ & $0.887$  & \multirow{2}{0.8cm}{fail} & $0.681$ & $0.727$ & $0.583$ & $0.665$ \\
          & $T$ & $208$ & $116$ &  & $36$ & $44$ & $26$ & $34$ \\
         \hline
    \end{tabular}
\caption{Performance comparison of the U-Net-based preconditioners. $\rho$ is the convergence factor and $T$ is the number of iterations required to reduce the residual norm by a factor of $10^6$. (a) summarizes the numerical results of the method on new examples equal in difficulty to the training set. Rows (b),(c) and (d) describe the performance of the same U-Nets on problems of a different size or with a different level of contrast in the slowness model. The ``single model'' column demonstrate the performance when the U-Net is trained using a single heterogeneous model. In the ``general model'' column the network is trained on multiple slowness models. The ``re-trained'' column describes the performance after mini-retraining.}
\label{table:convergence_factor_table}
\end{table}

\subsubsection{A stand-alone U-Net}
At the end of the training phase, the average relative error value for the stand-alone U-Net is $0.53$ (not shown). Due to the way the samples are produced and the limited network size, this value is consistent for both the training and the validation sets, i.e., there is no over-fitting. The preconditioner test of the U-Net is described in the top left plot in Fig. \ref{fig:retrain_results}. The preconditioner $M_{\sf VU}$ achieves the convergence within $46$ iterations and a convergence factor of $0.739$. The preconditioner $M_{\sf JU}$ achieves a convergence factor of $0.856$. When examining the network's ability to generalize to larger sizes than the training size, $M_{\sf JU}$ does not improve on the V-cycle performance, while $M_{\sf VU}$ reaches a convergence factor of $0.898$ for $256\times 256$-size and $0.972$ for $512\times 512$-size, compared to the V-cycle achieving $0.974$ and $0.989$, respectively. Next, we show that with little effort, the U-Nets can do better than that.

\subsubsection{An encoder-solver pair of networks}\label{sec:encoder-solver}
Taking advantage of the fact that the model does not change throughout the preconditioner test, we now use the encoder-solver framework, where another neural network pre-process the slowness model and generate context vectors for the solver U-Net. We train both networks as a single component, and in the preconditioner test we run the encoder network once and use the same feature vectors we have received from it throughout the iterations, only executing the solver network many times.

The average relative error of the training and validation sets at the end of the training phase is $0.49$. In the preconditioner test phase, the preconditioner $M_{\sf VU}$ reached convergence within $30$ iterations, and a convergence factor of $0.628$, while $M_{\sf JU}$ achieved $0.704$. Plot (d) of Fig. \ref{fig:retrain_results} demonstrates the results of the preconditioner test on $128\times 128$-grid inputs, consistent with the training set. The other two plots in (e) and (f) describe the generalization capability for larger-sized inputs and frequencies.

Beyond the improvement in the preconditioner test of $M_{\sf VU}$ and especially the more cost-effective $M_{\sf JU}$, the ability of the preconditioner to generalize on other sizes also improved significantly with the addition of the encoder network. Problems of $256 \times 256$-size converge within 62 iterations, compared to 391 using V-cycle only. For the $512 \times 512$ grids (with $\omega=80\pi$), the amount of iterations is reduced significantly---more than 600 iterations are required with $M_{\sf V}$ compared to 141 with $M_{\sf VU}$.

\subsubsection{Mini-retraining towards a single model network at solve time}
Due to the significant gap in the performance when generalizing over the slowness model (as opposed to knowing the model), we now add a mini-training phase to tune the network weights to a specific model. The experiment consists of three phases: training to learn the solver for different slowness models, mini-retraining to tune the network to a specific model, and the preconditioner test.

At this point, there is a trade-off between the cost of the retraining process and the success of the preconditioner test. The more resources we invest in the retraining phase, the closer the result will be to a single model result in Table \ref{table:convergence_factor_table}. We present the enhancement of a re-training phase that includes 30 epochs on 100 different $( \bfr,\bfe )$ pairs produced by a single slowness model $\kappa^2$. The re-training for problems of $256 \times 256$-size, was done using 30 epochs on 200 pairs while for problems from size of $512 \times 512$ we performed a re-training that includes 30 epochs on a set of 300 samples. These amounts of data samples and epochs are just a fraction of the corresponding amounts of a full neural network training. Still, the retraining is quite expensive and suits the case of multiple right-hand-sides only.
To reach the single model performance, a large set is required as re-training on a too small set may over-fit the input pairs and not just the model. The dotted lines in Fig. \ref{fig:retrain_results} graphs describe the effect of re-training on the preconditioner test.

In plot (d) of Fig. \ref{fig:retrain_results}, the encoder-solver network reached a good network model in its training, and hence performs well, and is only slightly improved by retaining, saving only 5-6 iterations in each case (20\%).
In contrast, when re-training for larger problems, the network learns both the model and the form of solution to the new problem size. The dotted lines in plots (e) and (f), show a significant improvement of the mini re-training process. For $256 \times 256$-size, the $M_{\sf VU}$ preconditioner converges within $46$ iterations compared to $62$ before the re-training. The $M_{\sf JU}$ preconditioner converges within $67$ iterations instead of $92$. For $512 \times 512$ problems, the improvement is even more significant, as we got convergence in 75 iterations compared to 141 iterations (almost twice) without the re-training process for learning the model and adjusting to $512 \times 512$-size. Note that SL preconditioner for $512\times 512$ converges in more than $600$ iterations, and $M_{\sf JU}$ also achieved a significant improvement following the re-training.

Plots (a), (b) and (c) in Fig. \ref{fig:retrain_results} show the results of the same experiment with a stand-alone U-Net. We see that the improvement is even more significant than in the encoder-solver case, probably because the starting point of the retraining phase is worse in the stand-alone experiment than in the encoder-solver case, and hence there is more room for improvement.

\subsubsection{The influence of the Krylov subspace and block sizes}
\begin{figure}
\centering
    \begin{tabular}{cc}
        \scriptsize{(a)}&\scriptsize{(b)}\\
        \includegraphics[width=0.47\textwidth]{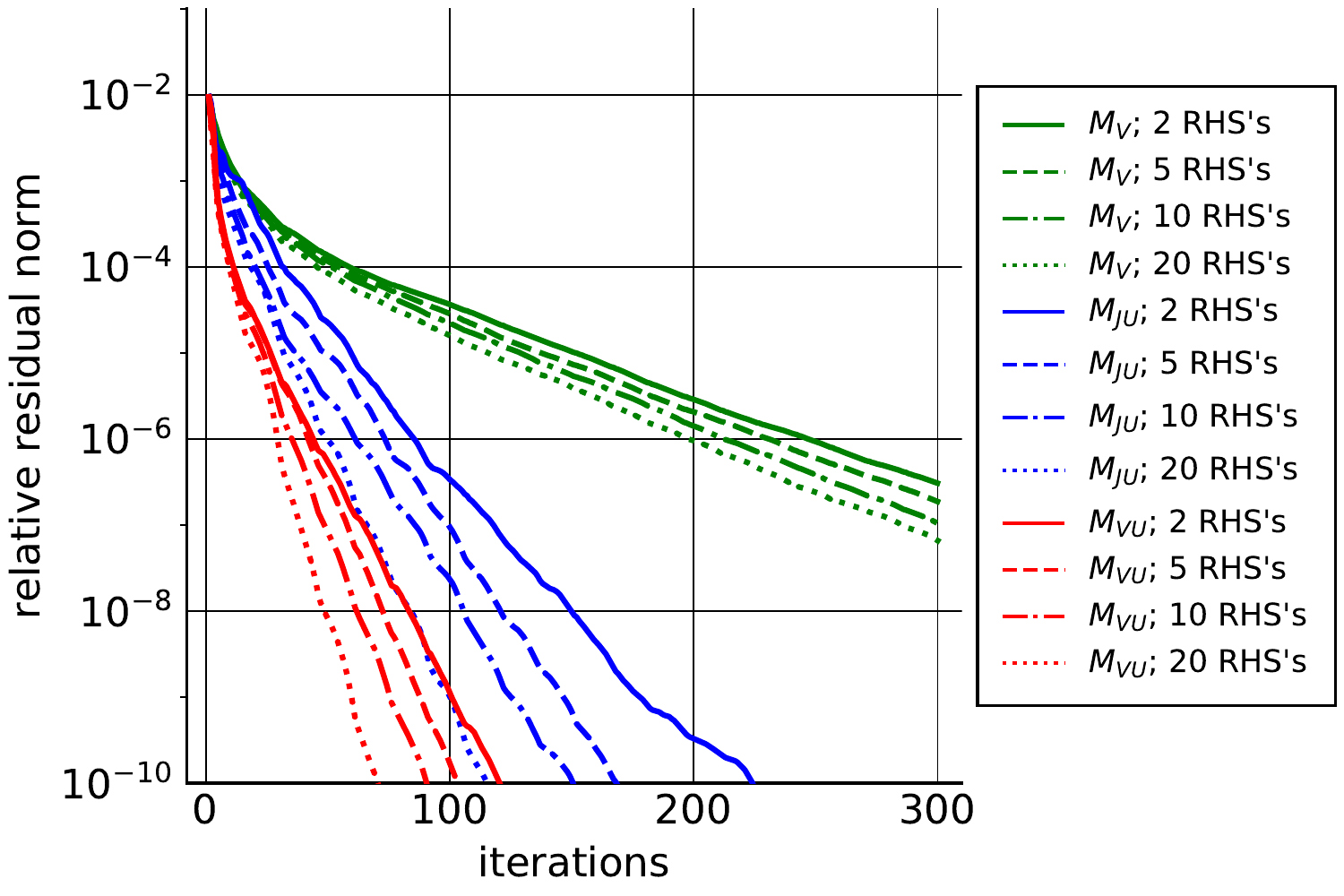}& \includegraphics[width=0.47\textwidth]{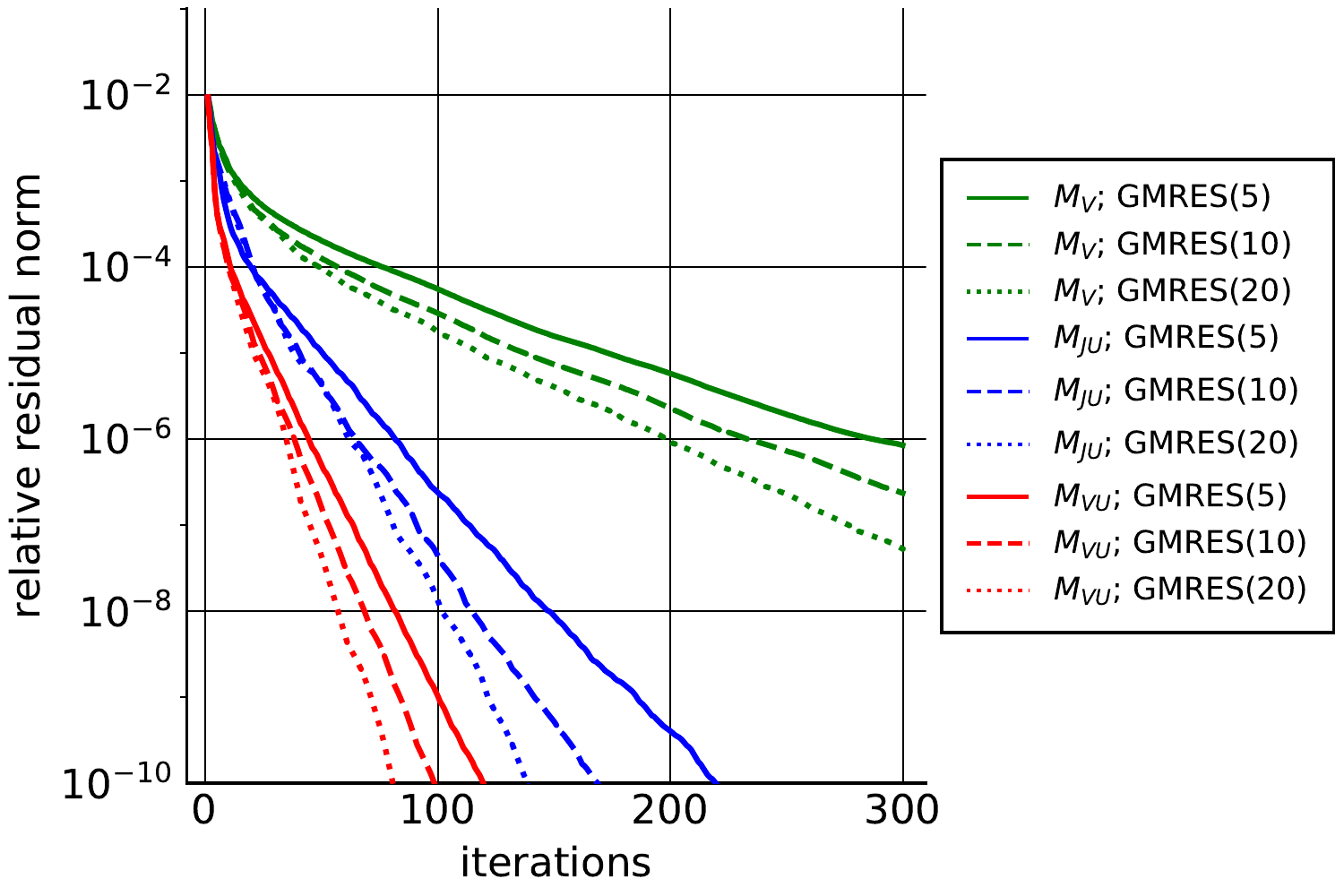}
    \end{tabular}
    \caption{The influence of the Krylov subspace and block sizes. On the left we see the influence of the block size (convergence for a number of right-hand-sides solved together by block FGMRES(10)). On the right we see the influence of the subspace size on the convergence factor. 
    }
    \label{fig:gmres_rhss}
\end{figure}

In this set of experiments we examine the influence of both the subspace size and number of RHSs solved together (block size) that the FGMRES method uses. The experiment corresponds to the $256\times256$ experiment shown in Fig \ref{fig:retrain_results} (e). Fig. \ref{fig:gmres_rhss} summerizes the results. It is obvious that whether by multiple RHSs or by subspace size, increasing the effective Krylov subspace can improve the performance of the solver significantly. That is in contrast to the SL V-cycle that is not improved by the same magnitude. This shows that there is more potential to the U-Net based preconditioners than shown earlier, since in addition, we typically did not reach the training error drop in the preconditioner test, indicating that FGMRES stabilizes the solution on the one hand, but also slows the preconditioners due to high residuals.

\subsubsection{GPU computation time}
In this set of experiments, we compare the running and convergence times of the different preconditioners. To have a fair comparison it is necessary to run the various preconditioners on the same processor unit. Since the U-Net is tailored to run efficiently on a GPU, we chose to have all the code running on a GPU as well, and compare the running times on this platform, which is quite common these days. In addition to GPU implementation of the multigrid cycle  described in Sec. \ref{sec:MG_GPU}, we also implemented the FGMRES method using the implementation of linear layers from the Julia Flux library. This function was used both in the external loop that compares performance of the preconditioners and in the coarse-grid solver of the V-cycle method. Since multiplication of 64-bit matrices and vectors on our GPU is significantly slower than those in 32-bit, we performed all the computations and timings in 32-bit. In our experiments, the convergence rate was not affected by the transition from CPU to GPU or from 64 to 32 bit, except to the minimum value of the residual norm which is limited to $10^{-8}$ only in 32-bit. This behavior is reasonable.

The FGMRES function that we use here supports block structure in the most naive way, which is less efficient than the full block-FGMRES used earlier in the CPU based implementation. That is because the block-FGMRES method includes QR decompositions which are highly inefficient on a GPU and are time consuming. As a result, as we increase the block size (number of right hand sides) the efficiency slightly drops but the computations better utilize the GPU resources. In our experiments we measured the times of running the algorithms on blocks of 10 RHSs. In addition, timings were averaged across 10 runs. The 32-bit GPU implementation was roughly 3-4 times faster than the CPU-based multicore implementation. Our experiments were conducted using an NVIDIA RTX 3090 GPU operated by the Julia Flux framework.

\begin{figure}
\centering
    \begin{tabular}{c c c}
        \scriptsize{(a)} & \scriptsize{(b)} & \scriptsize{(c)}\\
        \includegraphics[scale=0.35]{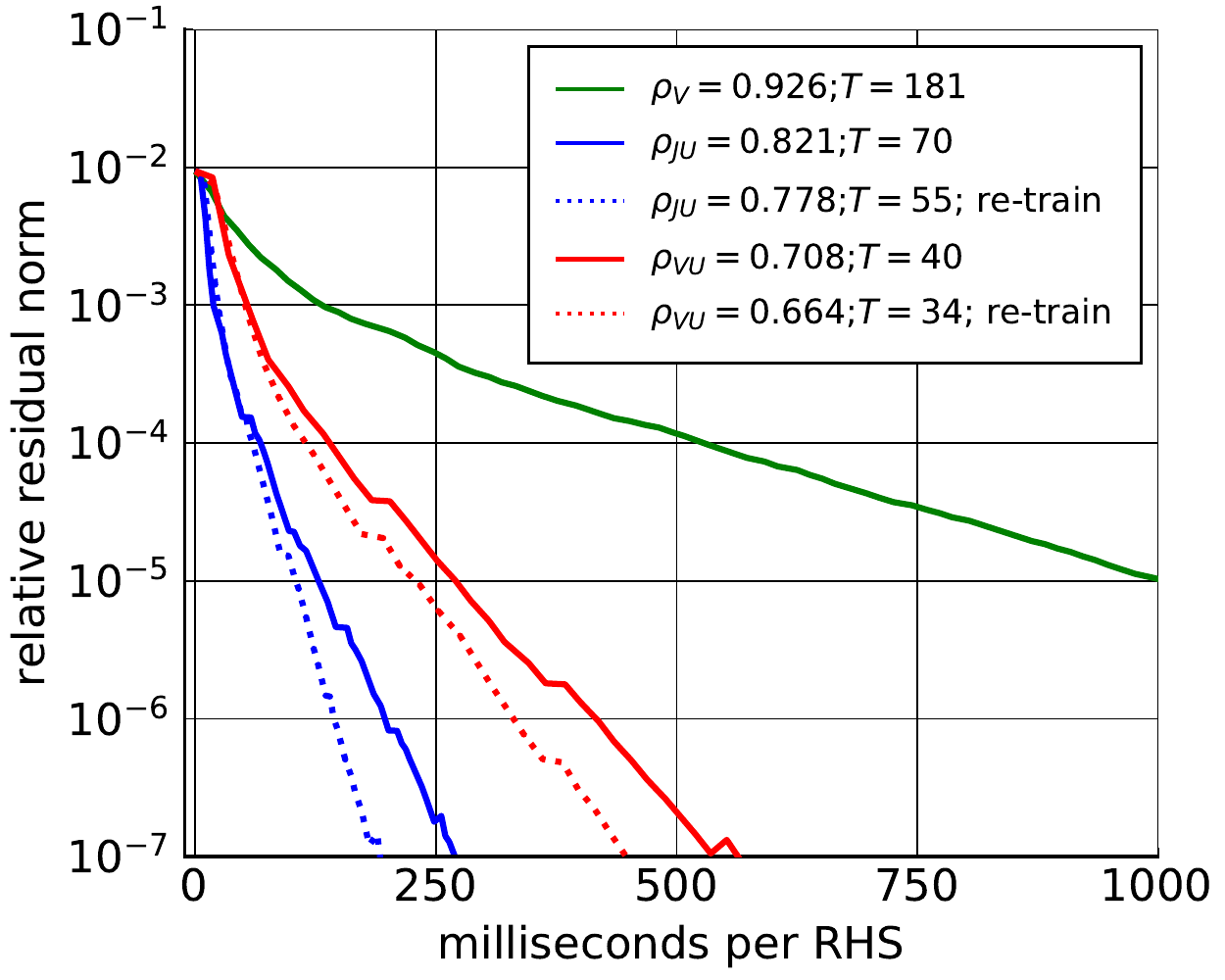} & \includegraphics[scale=0.35]{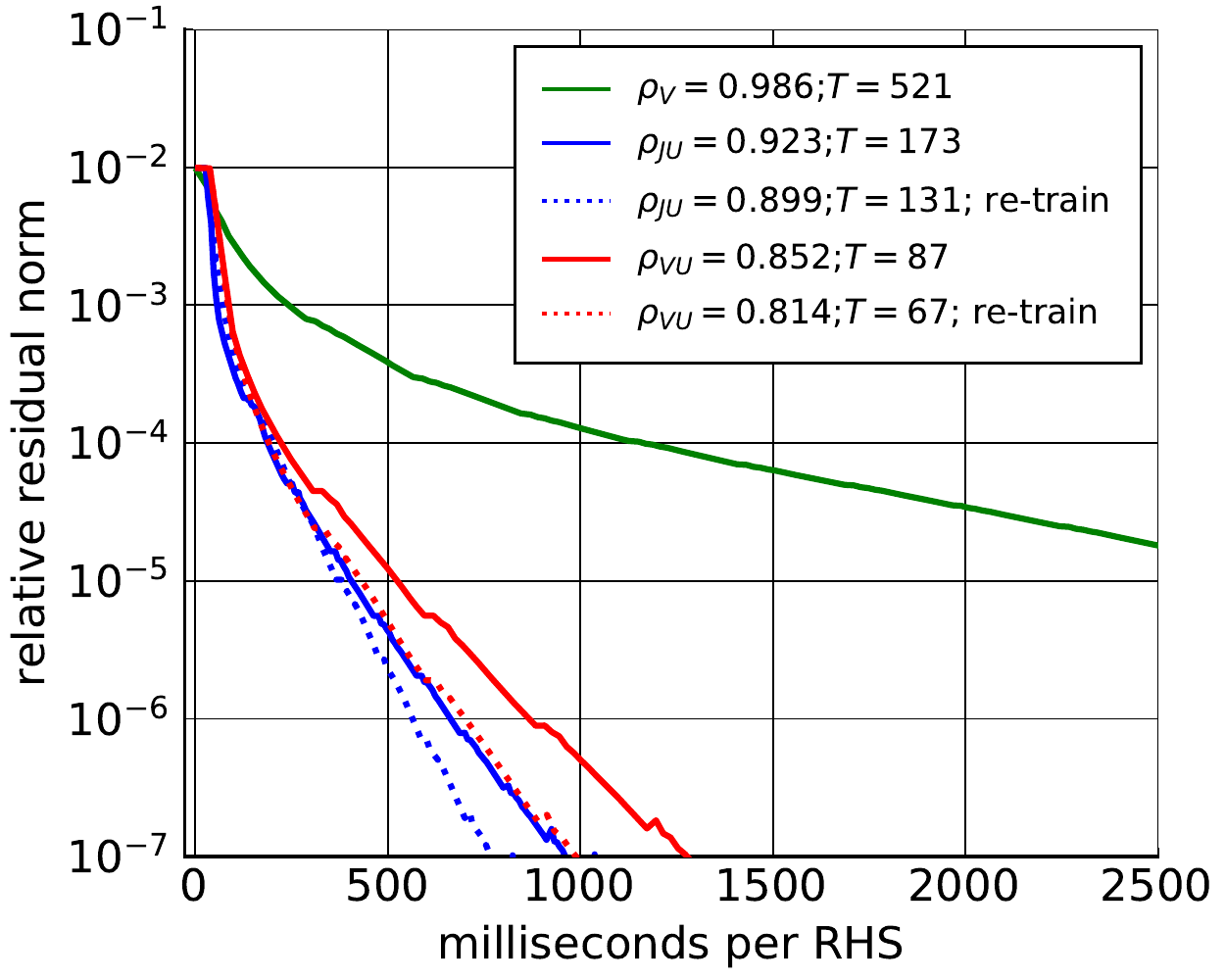} &
        \includegraphics[scale=0.35]{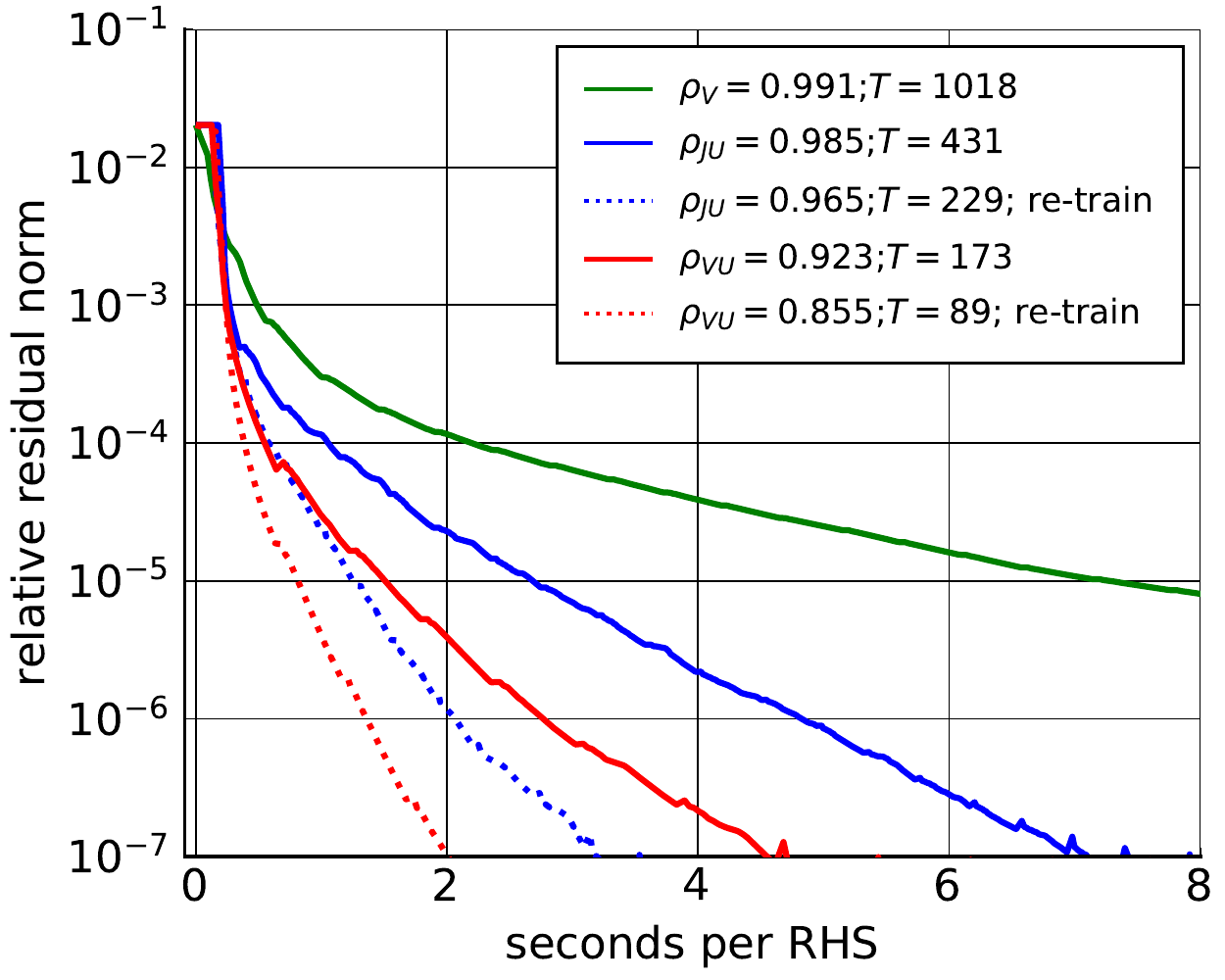}
    \end{tabular}
    \caption{Convergence history relative to times. Axis-x indicates the average time per RHS (in seconds / milliseconds) it took to reach each relative residual norm value. These experiments were measured on a GPU on 32-bit data describing problems of size (a) $128\times 128$, (b) $256 \times 256$ and (c) $512 \times 512$. 
    }
    \label{fig:time_results}
\end{figure}

Fig. \ref{fig:time_results} describes the residual convergence relative to the solution time. These runs correspond to the ones presented in the second row of Fig. \ref{fig:retrain_results} in terms of the problem setup (grid size, slowness models) for the various preconditioners, including those based on the re-train process. 
In addition, for each preconditioner we indicate the convergence factor (Eq. \eqref{eq:convergenceFact}), and the amount of iterations for the convergence.

In these experiments, the $M_{\sf JU}$ preconditioner for $128 \times 128$-grids as well as for $256 \times 256$-grids achieves better results than the $M_{\sf VU}$ preconditioner even though it performs a larger number of iterations, since the runtimes of each iteration are significantly smaller. Both, $M_{\sf JU}$ and $M_{\sf VU}$, achieve better results than the V-cycle, even for $512 \times 512$-grids and even without the re-training. For $128 \times 128$-size problems, the $M_{\sf JU}$ preconditioner after 250 milliseconds (per RHS) and the $M_{\sf VU}$ preconditioner after 600 milliseconds, reduced the residual norm by a factor of $10^6$, while the V-cycle-preconditioner reached the same result after more than 2.5 seconds. For $256 \times 256$-grids the two proposed preconditioners, $M_{\sf JU}$ and $M_{\sf VU}$,
convergence in less than 1.5 seconds while V-cycle converges after almost 10 seconds. For larger grids with size of $512 \times 512$, $M_{\sf JU}$ converged within 4 seconds, $M_{\sf VU}$ converged within 2.8 seconds, and V-cycle needed more than 25 seconds.

\subsubsection{Out of distribution models}

Having examined the ability of the U-Nets to solve problems produced using the CIFAR-10 training set of slowness models, in this experiment we present the ability of the CNN to adapt and generalize for realistic geophysical models that were not included in the training set. This aspect is important for the usability of the pre-trained CNN for new slowness models which are not available at training time. We do note that if one has an application in mind, it is better to use training data that are representative of that application. Here, we demonstrate a case where a new slowness model is out of the training distribution, as may inevitably happen in real life.

\begin{figure}
\centering
    \begin{tabular}{c c c}
        \scriptsize{(a)} & \scriptsize{(b)} & \scriptsize{(c)}\\
        \includegraphics[scale=0.38]{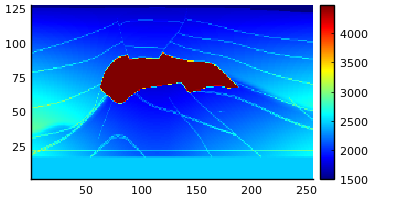} & \includegraphics[scale=0.38]{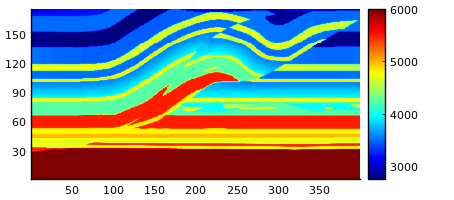} &
        \includegraphics[scale=0.38]{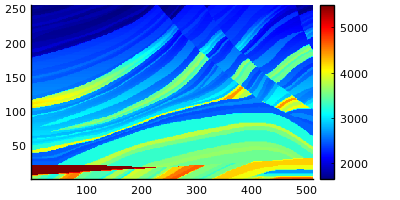}\\
        \scriptsize{(d)} & \scriptsize{(e)} & \scriptsize{(f)}\\
        \includegraphics[scale=0.35]{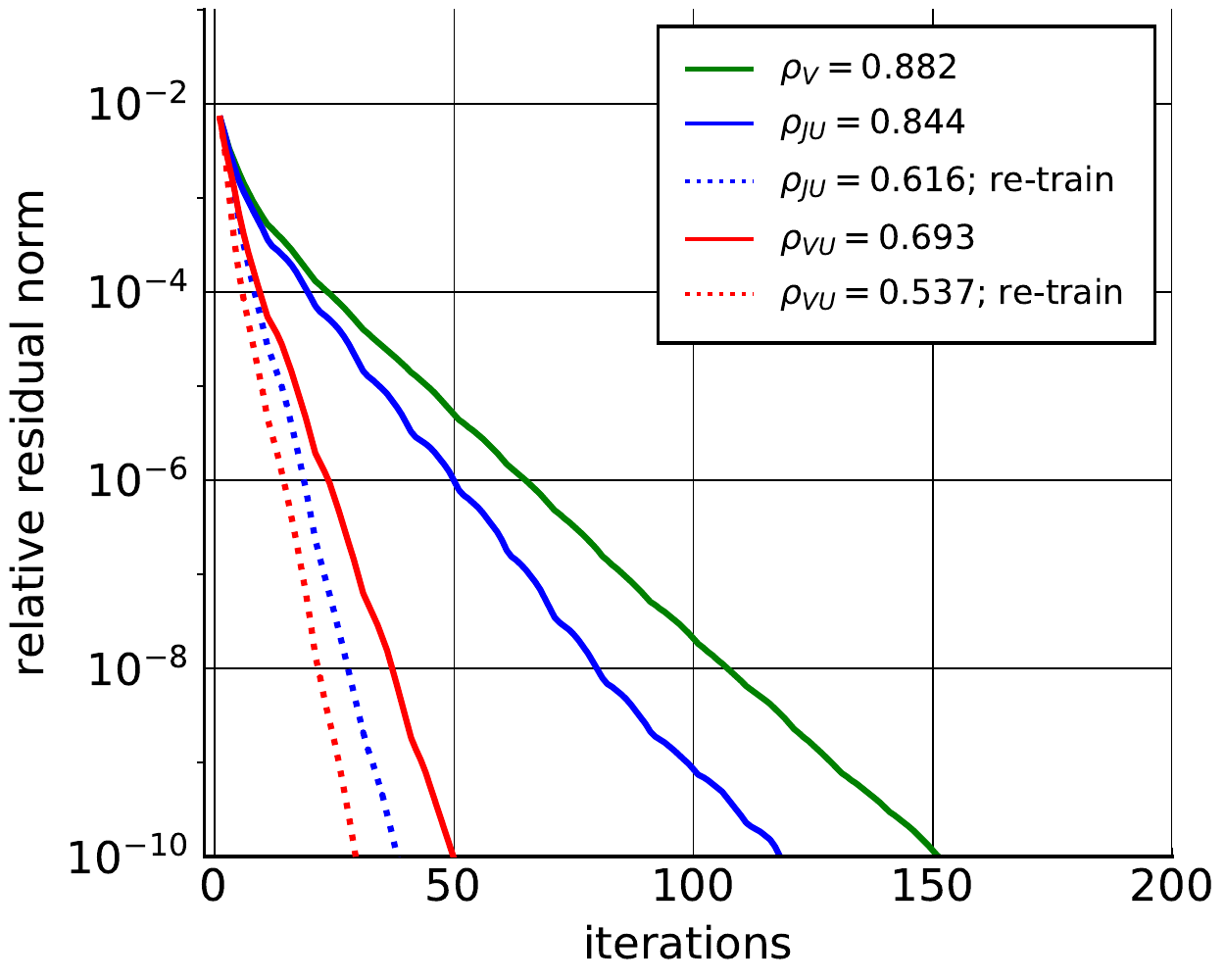} & \includegraphics[scale=0.35]{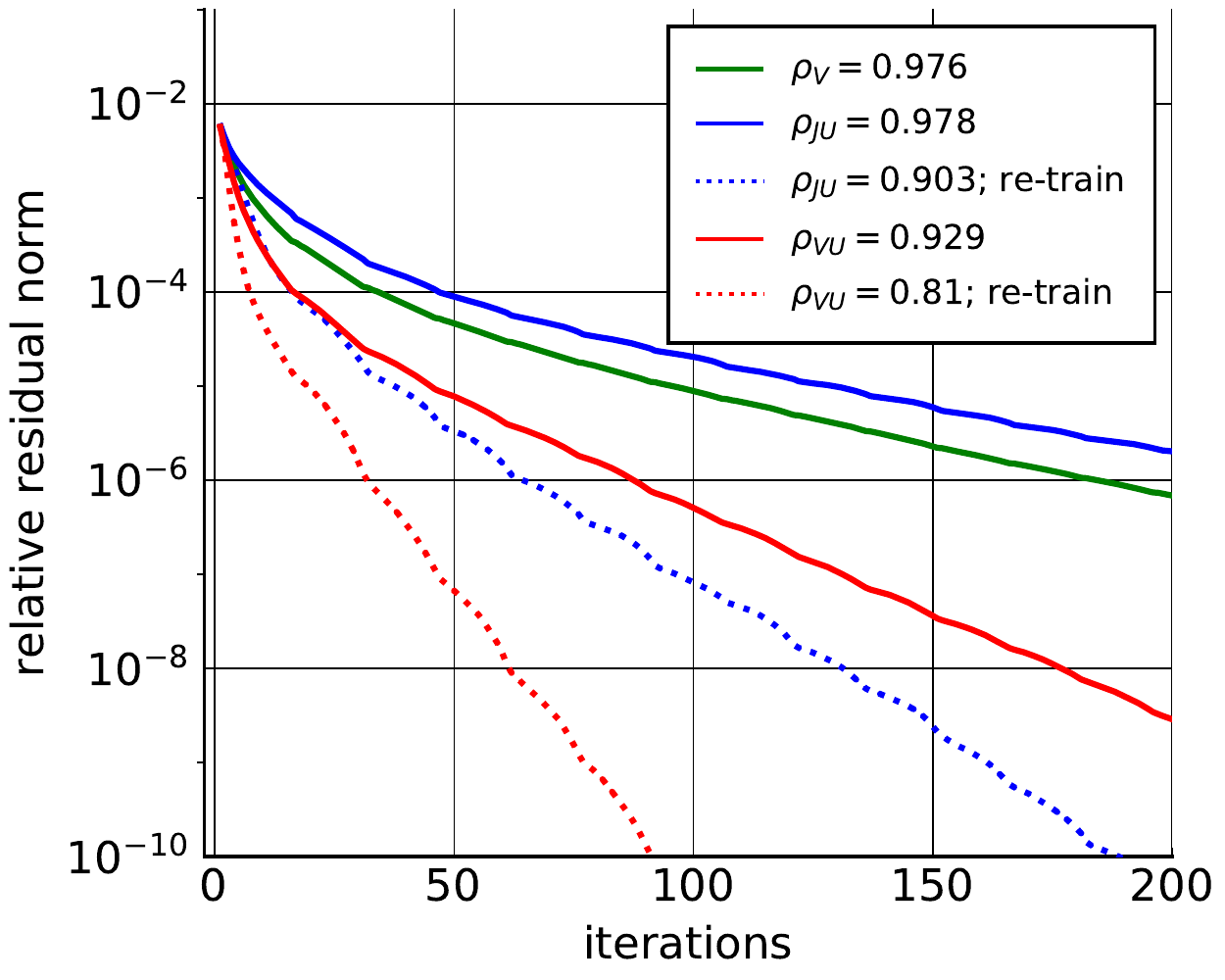} &
        \includegraphics[scale=0.35]{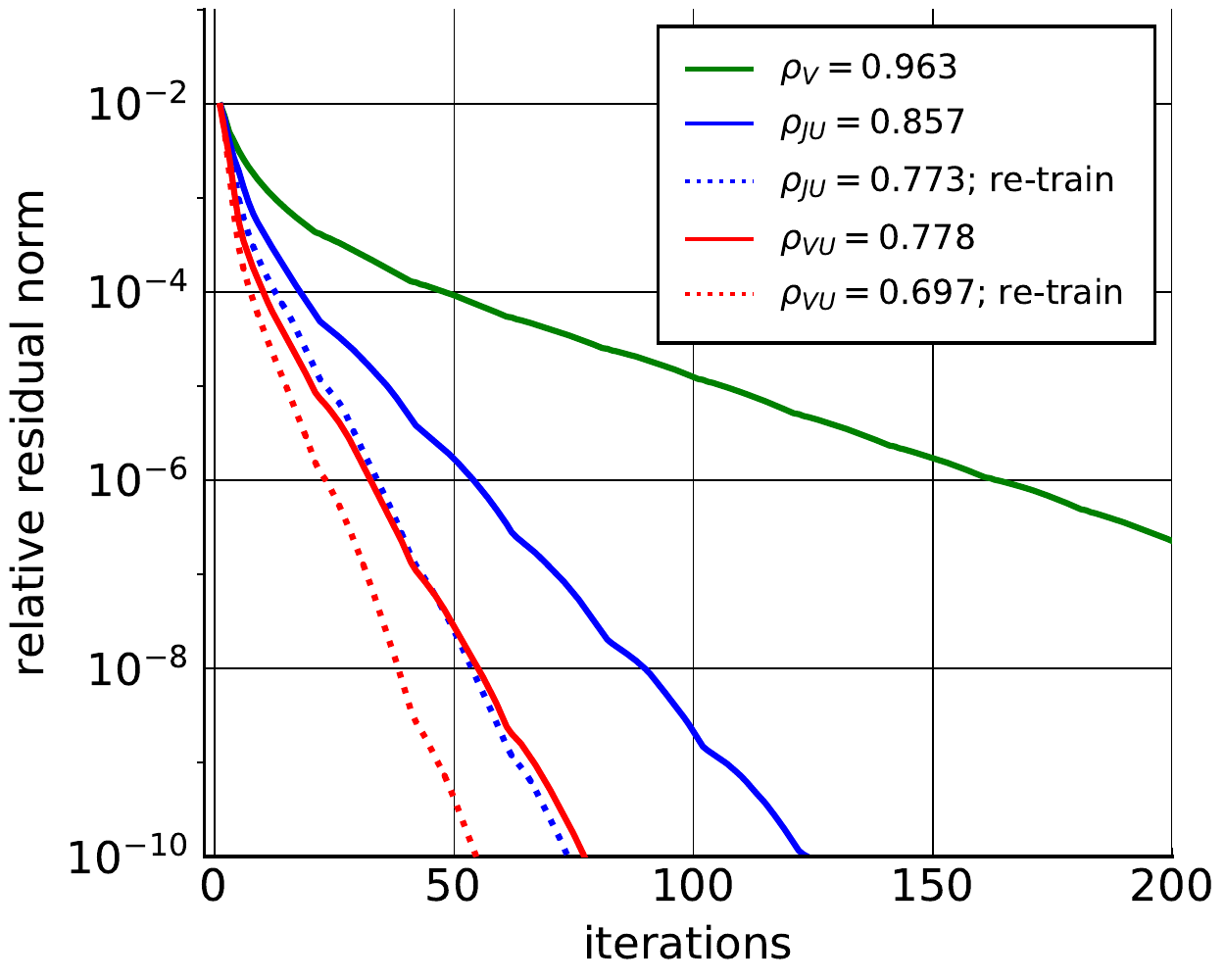}
    \end{tabular}
    \caption{Performance for geophysical models. First row: (a) SEG/EAGE salt-dome model ($128 \times 256$), (b) The Overthrust velocity model ($176 \times 400$) and (c) the Marmousi velocity model ($256 \times 512$). Second row: the convergence history of block-FGMRES
    on 10 Helmholtz problems produced by the respective slowness model from the row above. The dotted lines describe the results of a 30 epochs re-training process for the specific model on a training set of 300 RHSs.}
    \label{fig:ood_results}
\end{figure}

Fig. \ref{fig:ood_results} describes the comparison between the performance of the encoder-solver U-Nets and V-cycle  for problems produced by 3 different geophysical slowness models that are commonly used in the literature---the SEG salt and Overthrust models \cite{aminzadeh19973}, and the Marmousi model \cite{brougois1990marmousi}. These models include many sharp edges and are significantly different than the rather smooth models that we produced to train the CNN. We again compare the convergence history of the residual norm with respect to the iterations of the block-FGMRES for each of the preconditioners. As before, the combination of U-Net and V-cycle ($M_{\sf VU}$) is more effective than the V-cycle alone even without the re-training phase, i.e. a CNN trained on other types of models of size $ 128 \times 128$, effective for different geophysical models and problem sizes not seen during training.

For each of the models we did a re-training that included 30 epochs on a set of 300 pairs $(\bfr,\bfe)$ produced using the specific slowness model, and marked their performance with dotted lines in the plots. First, for the three models the re-training process has made the $M_{\sf JU}$ preconditioner more effective than the V-cycle. In particular, for the first two models, the SEG/EAGE salt-dome model and the Overthrust velocity model, re-training saved about 100 iterations, depending on the desired threshold of the residual norm. The cost of the retraining: 30 (epochs) $\cdot$ 300 (train-set size) $\cdot$ 3 (U-Net operation $+$ gradient calculation + the encoder's gradient) amounts to about 27,000 U-Net applications. Divided by the gain of 100 iterations per RHS, the retraining \emph{in these settings} will be beneficial if we need to solve about 270 RHSs. Such a number is small compared to the number of RHSs that need to be solved as part of a solution of an inverse problem with many sources. Note that the slowness models here are significantly different than the ones used for training, and we expect a lighter retraining phase if the grid-size and training set distribution better match the target problems.

\section{Conclusions}
In this work we presented a new approach for preconditioning the heterogenouos Helmholtz equation. Our approach is based on a combination of CNN and multigrid components, using a U-Net CNN as a coarse grid solver, or as a deflation operator together with SL V-cycles. The idea is that the CNN is trained to reduce the error, and the smoothing or V-cycles are responsible for stabilizing the solution process in a Krylov solver.

Since we target a preconditioner, our CNN is trained to generalize over the right-hand-sides (residuals), and over the slowness models.
To help with the generalization, we propose to use an encoder-solver framework, where two CNNs are used. One prepares context vectors and generalizes over the slowness, and the other receives the context vectors and generalize over the right-hand-sides. This encoder-solver framework indeed upgrades the stand-alone U-Net in all cases. In addition, inspired by PINNs, we offer a mini-retraining phase over the slowness model at solve time, since the method is more effective when the slowness model is known. This consistently shows yet another upgrade to the performance. We show that while our U-Net may require more FLOPs than traditional methods, it can be applied efficiently on GPU hardware, and yield favorable running times. Furthermore, for many cases CNN models can be significantly compressed by a variety of ways without harming their performance, making the preconditioner more beneficial.

In this work we used the U-Net architecture ``as is''. However, better architectures can probably be defined specifically for the Helmholtz equation. Furthermore, since problems typically include a more specific right-hand-side (e.g., point source or a plane wave), the first iteration can be performed with a different network that is trained for such cases, and our framework can then follow to improve the solution. Another interesting direction can be to train iterative U-Nets directly (like Recurrent Neural Networks), removing the need of FGMRES to guide the solution.

Lastly, in this paper we demonstrated 2D problems only, so that the number of waves in the domain is high and realistic. Performing the same in 3D will require the training of 3D networks in the same resolution, which is beyond the abilities of our hardware at this point. Small 3D problems are solved by V-cycles quite efficiently, hence are less relevant. However, it is indeed relevant for HPC centers solving large 3D problems, because the training is done offline.

\bibliographystyle{siamplain}
\bibliography{Research_Paper}
\end{document}